\documentclass[review,12pt]{elsarticle}

\usepackage{amssymb}
\usepackage{amsthm}
\usepackage{framed} 
\usepackage{amsmath,color}
\usepackage{mathrsfs}
\usepackage{graphicx}
\usepackage{epstopdf}
\usepackage{float}
\usepackage{caption}
\usepackage{subcaption}
\usepackage{bm}
\usepackage{bbm}
\usepackage{mathrsfs}
\usepackage{cleveref}
\usepackage{soul}
\usepackage{accents}
\usepackage{color,soul} 
\usepackage{color} 
\usepackage{bm}
\usepackage{multirow} 
\biboptions{sort&compress}
\soulregister\citep7 
\soulregister\citet7 
\soulregister\citealp7 
\newsavebox{\measurebox} 
\usepackage{titlesec} 
\usepackage[T1]{fontenc} 
\usepackage{lmodern} 
\input{glyphtounicode}  

\newcommand{\tm}{\textrm} 

\def\onedot{$\mathsurround0pt\ldotp$}
\def\cddot{
  \mathbin{\vcenter{\baselineskip.67ex
    \hbox{\onedot}\hbox{\onedot}}%
  }}%
\def\cdddot#1{
  \mathbin{\vcenter{\baselineskip.67ex
    \hbox{\onedot}\hbox{\onedot}\hbox{\onedot}%
  }}%
}

\journal{Theoretical and Applied Fracture Mechanics}

\makeatletter
\def\@author#1{\g@addto@macro\elsauthors{\normalsize%
    \def\baselinestretch{1}%
    \upshape\authorsep#1\unskip\textsuperscript{%
      \ifx\@fnmark\@empty\else\unskip\sep\@fnmark\let\sep=,\fi
      \ifx\@corref\@empty\else\unskip\sep\@corref\let\sep=,\fi
      }%
    \def\authorsep{\unskip,\space}%
    \global\let\@fnmark\@empty
    \global\let\@corref\@empty  
    \global\let\sep\@empty}%
    \@eadauthor={#1}
}
\makeatother

\titleformat{\paragraph}
{\normalfont\normalsize\itshape}{\theparagraph}{1em}{}
\titlespacing*{\paragraph}
{0pt}{3.25ex plus 1ex minus .2ex}{1.5ex plus .2ex}

\begin{document}

\begin{frontmatter}



\title{Phase field fracture modelling using quasi-Newton methods and a new adaptive step scheme}


\author{Philip K. Kristensen \fnref{DTU}}


\author{Emilio Mart\'{\i}nez-Pa\~neda\corref{cor1}\fnref{IC}}
\ead{e.martinez-paneda@imperial.ac.uk}

\address[DTU]{Department of Mechanical Engineering, Technical University of Denmark, DK-2800 Kgs. Lyngby, Denmark}


\address[IC]{Department of Civil and Environmental Engineering, Imperial College London, London SW7 2AZ, UK}

\cortext[cor1]{Corresponding author.}

\begin{abstract}
We investigate the potential of quasi-Newton methods in facilitating convergence of monolithic solution schemes for phase field fracture modelling. Several paradigmatic boundary value problems are addressed, spanning the fields of quasi-static fracture, fatigue damage and dynamic cracking. The finite element results obtained reveal the robustness of quasi-Newton monolithic schemes, with convergence readily attained under both stable and unstable cracking conditions. Moreover, since the solution method is unconditionally stable, very significant computational gains are observed relative to the widely used staggered solution schemes. In addition, a new adaptive time increment scheme is presented to further reduces the computational cost while allowing to accurately resolve sudden changes in material behavior, such as unstable crack growth. Computation times can be reduced by several orders of magnitude, with the number of load increments required by the corresponding staggered solution being up to 3000 times higher. Quasi-Newton monolithic solution schemes can be a key enabler for large scale phase field fracture simulations. Implications are particularly relevant for the emerging field of phase field fatigue, as results show that staggered cycle-by-cycle calculations are prohibitive in mid or high cycle fatigue. The finite element codes are available to download from www.empaneda.com/codes.
\end{abstract}


\end{frontmatter}



\section{Introduction}
\label{Sec:Introduction}

The phase field fracture method has emerged as a promising variational framework for modelling advanced cracking problems. Fracture can be revisited as an energy minimisation problem by solving for an auxiliary variable, the phase field parameter $\phi$  \cite{Francfort1998,Bourdin2000}. Consequently, complex fracture features, such as crack branching, crack initiation from arbitrary sites or coalescence of multiple cracks, are naturally captured in the original finite element mesh (see, e.g., \cite{Borden2012,Borden2016,McAuliffe2016}). Not surprisingly, the method is becoming increasingly popular and the number of applications has soared. Recent examples include hydrogen embrittlement \cite{CMAME2018,Duda2018}, fatigue damage \cite{Alessi2018c,Carrara2020}, cracking of lithium-ion batteries \cite{Miehe2015,Zhao2016}, rock fracture \citep{Zhou2018}, composites delamination \cite{Quintanas-Corominas2019,Quintanas-Corominas2020}, and fracture of functionally graded materials \citep{CPB2019}, among other; see \cite{Wu2020} for a review.  \\

A great deal of attention has been devoted to the development of efficient schemes for solving the coupled deformation-fracture problem. The total potential energy functional, including the contributions from the bulk and fracture energies, is non-convex with respect to the primary kinematic variables, the displacement field $\bm{u}$ and the phase field $\phi$. Due to this non-convexity, the Jacobian matrix in Newton's method becomes indefinite, hindering convergence and robustness in monolithic solution schemes, where $\bm{u}$ and $\phi$ are solved simultaneously. Different numerical strategies have been adopted to overcome these drawbacks: error-oriented Newton methods \cite{Wick2017a}, \textit{ad hoc} line search algorithms \cite{Gerasimov2016,Singh2016,Borst2016a} and modified Newton methods \cite{Wick2017}. While promising results have been obtained, performance is very problem-dependent and the monolithic minimisation of the energy functional ``remains extremely challenging'' \cite{Wick2017}. Staggered solution schemes, based on alternating minimisation, enjoy a greater popularity \cite{Miehe2010,Miehe2010a,Hofacker2012,Ambati2015,Kristensen2020}. By fixing one primal kinematic variable, the total potential energy becomes convex with respect to the other primal kinematic variable. The method has proven to be very robust but computationally demanding. First, convergence of critical loading steps requires a significant amount of iterations \cite{Gerasimov2016}. In addition, the method is no longer unconditionally stable, requiring the use of very small load steps to effectively track the equilibrium solution \cite{Singh2016} or recursive iteration schemes \cite{Bourdin2000}.\\

In this work, we demonstrate that a robust and efficient numerical framework for phase field fracture analyses can be obtained by combining quasi-Newton methods and a monolithic solution scheme. There is a large literature devoted to the analysis of the robustness of quasi-Newton methods when dealing with non-convex minimization problems - see, e.g., \cite{Li2001,Li2001b,Lewis2013} and references therein. Very recently, Wu \textit{et al.} \cite{Wu2020a} showed the potential of quasi-Newton monolithic approaches in the context of the so-called unified phase field damage theory, a phase field regularisation of cohesive zone models (PF-CZM) \cite{Wu2017,Wu2018a}. We extend their analysis to the standard phase field fracture formulation and showcase the potential of the method in three problems of different nature: quasi-static fracture, phase field fatigue and dynamic fracture. In addition, we introduce a new adaptive time stepping criterion for phase field cracking. The results obtained reveal computation times that are up to 100 times smaller than those required to obtain the same result with the widely used staggered solution. These results back the earlier findings by Wu \textit{et al.} \cite{Wu2020a} in the context of quasi-static fracture and the PF-CZM model, emphasising the promise of monolithic quasi-Newton implementations for phase field fracture and fatigue modelling.\\

The remainder of this manuscript is organized as follows. The theoretical phase field formulation employed to model (quasi-static and dynamic) fracture is shown in Section \ref{Sec:Theory}. Details of the numerical implementation are given in Section \ref{Sec:Numerical}, including a comprehensive description of the Broyden-Fletcher-Goldfarb-Shanno (BFGS) algorithm employed. Representative numerical results are shown in Section \ref{Sec:Results}. First, paradigmatic examples in quasi-static phase field fracture are revisited with the new solution scheme. The analysis is then extended to the case of fatigue cracking. And finally, the potential of the method is showcased in dynamic fracture, where off-diagonal matrices have a larger relative weight. Concluding remarks are given in Section \ref{Sec:Conclusions}.

\section{The phase field fracture method}
\label{Sec:Theory}

\subsection{Phase field approximation of the fracture energy}

\begin{figure}[H] 
    \centering
    \includegraphics[scale=0.9]{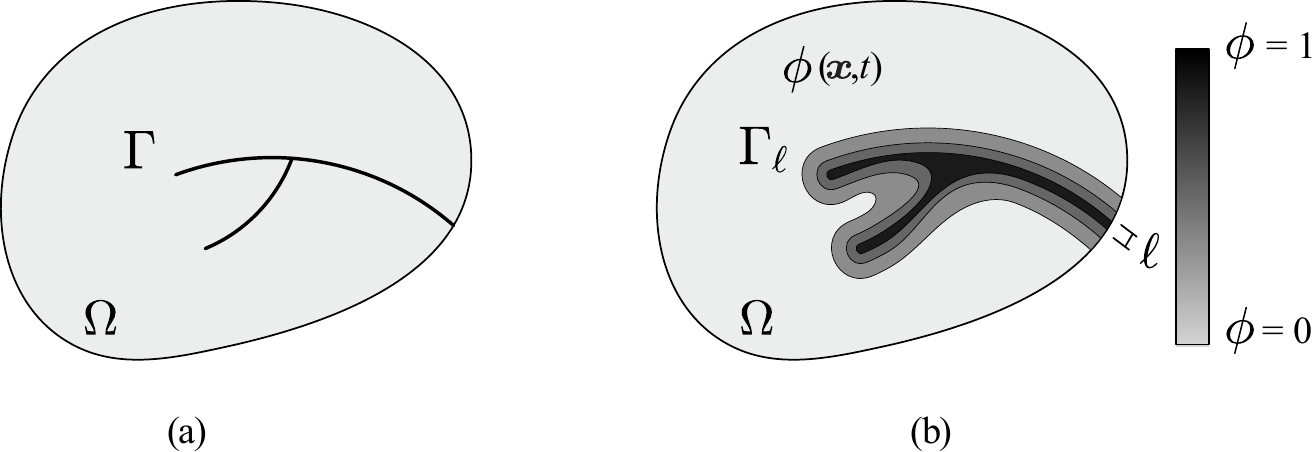}
    \caption{Schematic representation of a solid body with (a) internal discontinuity boundaries, and (b) a phase field approximation of the discrete discontinuities.}
    \label{fig:PhaseFieldFigIntro}
\end{figure}

Alan Arnold Griffith's energy-based analysis of cracks in 1920 is considered to be the birth of the field of fracture mechanics \cite{Griffith1920}. Consider a cracked solid with strain energy density $\psi (\bm{\varepsilon})$, which is a function of the strain tensor $\bm{\varepsilon}$. In the absence of external forces, the variation of the total energy $\Pi$ due to an incremental increase in the crack area d$A$ is given by
\begin{equation}
\frac{\text{d} \Pi}{\text{d} A} = \frac{\text{d} \psi (\bm{\varepsilon})}{\text{d} A} + \frac{\text{d} W_c}{\text{d} A}  = 0,
\end{equation}

\noindent where $W_c$ is the work required to create new surfaces. The last term is the so-called critical energy release rate $G_c=\text{d} W_c / \text{d} A$, a material property that characterises the fracture resistance. Griffith's energy balance can be formulated in a variational form as:
\begin{equation}\label{Eq:Pi}
\Pi = \int_\Omega \psi \left( \bm{\varepsilon} \right) \text{d} V + \int_\Gamma   G_c \, \text{d} \Gamma,
\end{equation}

\noindent with $\Gamma$ being the crack surface and $V$ denoting the volume of the solid, occupying an arbitrary domain $\Omega$. The crack surface is unknown, hindering minimization of (\ref{Eq:Pi}). However, an auxiliary variable, the phase field $\phi$, can be used to track the crack interface, see Fig. \ref{fig:PhaseFieldFigIntro}. The phase field $\phi$ is a damage-like variable that takes the values of 0 in an intact material point, and of 1 in a fully cracked material point.

Following continuum damage mechanics arguments, a degradation function $g \left( \phi \right) = \left( 1 - \phi \right)^2$ is defined that diminishes the stiffness of the material with evolving damage \cite{Bourdin2000}. Accordingly, the total potential energy functional can be formulated as
\begin{equation}\label{Eq:Piphi}
\Pi_\ell = \int_\Omega \left\{ \left( 1 - \phi \right)^2 \psi_0 \left( \bm{\varepsilon} \right) + G_c \left( \frac{1}{2 \ell} \phi^2  + \frac{\ell}{2} |\nabla \phi|^2 \right) \right\} \text{d} V
\end{equation}
\noindent where $\ell$ is a length scale parameter that governs the size of the fracture process zone and $\psi_0$ denotes the elastic strain energy of the undamaged solid. The work required to create a cracked surface, $\Gamma$, is now expressed as a volume integral, making the problem computationally tractable. As shown by $\Gamma$-convergence, the regularized functional $\Pi_\ell$ approaches the functional of the discrete crack problem $\Pi$ for $\ell \to 0$ \cite{Bellettini1994,Chambolle2004}. The choice (\ref{Eq:Piphi}) is based on the work by Bourdin \textit{et al.} \cite{Bourdin2000} and the earlier regularization by Ambrosio and Tortorelli of the Mumford-Shah problem in image processing \cite{Ambrosio1991}. This surface regularization is commonly referred to as the AT2 model. See Ref. \cite{Mandal2019} for other choices and a detailed numerical comparison in the context of phase field fracture. Considering the earlier work by Wu \textit{et al.} \cite{Wu2020a}, the superior performance of monolithic quasi-Newton solution strategies is therefore demonstrated for both the PF-CZM and AT2 regularizations; the analysis of the so-called AT1 model \cite{Pham2011} remains to be addressed.

\subsection{Governing balance equations of the coupled problem}
\label{2.4}
\subsubsection{Basic fields and boundary conditions}

\begin{figure}[H]
    \centering
    \includegraphics[scale=1.2]{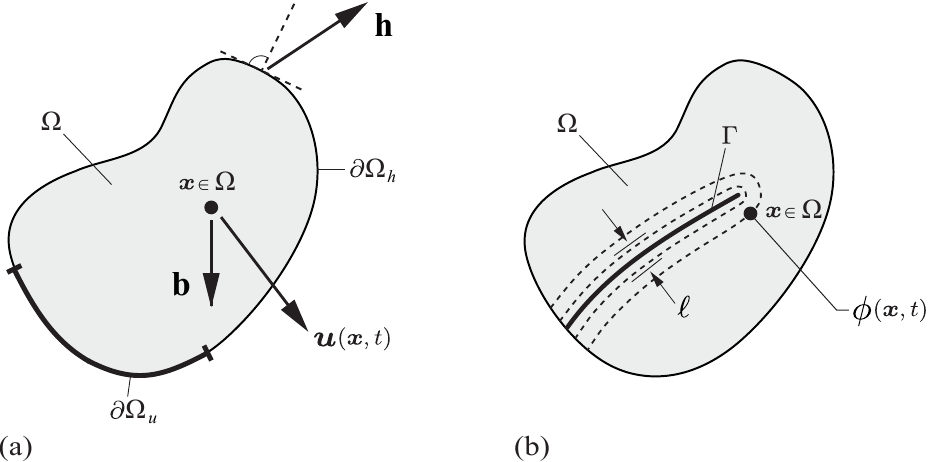}
    \caption{Two-field model of phase field fracture in deformable solids. The displacement field $\bm{u}$ is constrained by the Dirichlet- and Neumann-type boundary conditions \(\bm{u}=\bm{u}_{\text{D}}\,\) on \(\,\partial \Omega_{u}\) and \(\bm{\sigma} \cdot \bm{n}=\bm{h} \,\) on \(\, \partial \Omega_{h}\). (b) The crack phase field \(\phi\) is constrained by the Dirichlet- and Neumann-type boundary conditions \(\phi=1 \,\) on \(\, \Gamma\) and \(\nabla \phi \cdot \bm{n} =0\,\) on \(\, \partial \Omega\).}
    \label{fig:BC}
\end{figure}

We proceed now to formulate the governing equations for the displacement field \(\bm{u}\) and the phase field \(\phi\). With respect to \(\bm{u}\), the outer surface of the body is decomposed into a part \(\partial \Omega_{u}\), where the displacement is prescribed by Dirichlet-type boundary conditions
\begin{equation}
    \bm{u}(\bm{x},t)=\bm{u}_{\text{D}}(\bm{x},t) \hspace{5mm} at \hspace{3mm} \bm{x} \in \partial \Omega_{u}
    \centering
\end{equation}
and into a part \(\partial \Omega_{h}\), where the traction \(\bm{h}\) is prescribed by Neumann-type boundary conditions (see Fig. \ref{fig:BC}a). With respect to the fracture phase field, a cracked region can be prescribed through the initial condition
\begin{equation}
    \phi(\bm{x},t)=1 \hspace{5mm} at \hspace{3mm} \bm{x} \in \Gamma_{\text{D}}
    \centering
\end{equation}
\noindent where $\Gamma_{\text{D}}$ is a possible given sharp crack surface inside the solid $\Omega$ (see Fig. \ref{fig:BC}b). The crack phase field $\phi$ is considered to be driven by the displacement field $\bm{u}$ of the solid. Consequently, no prescribed external loading is considered corresponding to the crack phase field $\phi$.

As both quasi-static and dynamic fracture will be considered in this work, we define the kinetic energy of the solid as:
\begin{equation} \label{eq:Ekin}
\Psi^k \left( \dot{\bm{u}} \right) = \frac{1}{2} \int_\Omega \rho \, \dot{\bm{u}}\cdot\dot{\bm{u}} \, \text{d}V,
\end{equation}
where $\rho$ is the material density and $\dot{\bm{u}}=\partial\bm{u}/\partial t$. 
\subsubsection{Coupled balances}
With the kinetic and potential energies defined, along with the boundary conditions of the system, the Lagrangian for the regularised fracture problem is given by
\begin{align}\label{eq:Lagrangian}
    L\left(\bm{u},\dot{\bm{u}},\phi\right) = \Psi^k \left( \dot{\bm{u}} \right)-\Pi_\ell\left(\bm{u},\phi\right) .
\end{align} 
By insertion of (\ref{Eq:Piphi}) and (\ref{eq:Ekin}) into (\ref{eq:Lagrangian}), the Lagrangian can be formulated as: 
\begin{align}
    L\left(\bm{u},\dot{\bm{u}},\phi\right) = \int_\Omega \left\{\frac{1}{2}\rho \dot{\bm{u}}\cdot\dot{\bm{u}} - \left( 1 - \phi \right)^2 \psi_0 \left( \bm{\varepsilon} \right) - G_c \left( \frac{1}{2 \ell} \phi^2  + \frac{\ell}{2} |\nabla \phi|^2 \right)\right\}\, \tm{d}V. 
\end{align}
The weak form can be readily obtained by taking the stationary of the Lagrangian functional $\delta L\left(\bm{u},\dot{\bm{u}},\phi\right) =0$, such that:
\begin{align}\label{eq:weak0}
        \int_{\Omega} & \left\{\rho \ddot{\bm{u}} \cdot \delta \bm{u}+ \bm{\sigma} \cddot \delta \bm{\varepsilon} -2(1-\phi)\delta \phi \psi_{0}(\bm{\varepsilon}) +
        G_{c}\left(\dfrac{1}{\ell}\phi \delta \phi
        + \ell\nabla \phi \cdot \nabla \delta \phi \right) \right\}  \, \mathrm{d}V \nonumber \\ & -
        \int_{\Omega} \bm{b} \cdot \delta \bm{u}\, \mathrm{d}V - \int_{\partial \Omega_{h}} \bm{h} \cdot \delta \bm{u}\, \mathrm{d}S = 0.
\end{align}
\noindent Here, $\bm{b}$ is a prescribed body force field per unit volume and the Cauchy stress tensor $\bm{\sigma}$ is given in terms of the stress tensor of the undamaged solid $\bm{\sigma}_0$ and the degradation function $g (\phi)$ as:
\begin{equation}
    \bm{\sigma}=g(\phi)\bm{\sigma_{0}}=(1-\phi)^{2} \bm{\sigma_{0}} =  \bm{C_{0}} \cddot \bm{\varepsilon},
    \centering
\end{equation}
\noindent with $\bm{C_{0}}$ being the linear elastic stiffness matrix. By application of the Gauss divergence theorem and considering that (\ref{eq:weak0}) must hold for any arbitrary permissible variations $\delta \bm{u}$, $\delta\dot{\bm{u}}$ and $\delta\phi$, we arrive at the balance equations:
\begin{equation} \label{eq:strong}
    \begin{split}
        \nabla \bm{\sigma} +\bm{b} &= \rho \ddot{\bm{u}} \\[1mm]
        G_{c} \left(\dfrac{1}{\ell} \phi - \ell \Delta \phi \right) - 2(1-\phi) \psi_{0}(\bm{\varepsilon}) &= 0
    \end{split}
    \centering
\end{equation}
\noindent where \(\Delta \phi\) refers to the Laplacian of the phase field. 
The strong form balance eqautions are subject to the Neumann-type boundary conditions
\begin{equation} \label{eq:strongBC}
    \bm{\sigma} \cdot \bm{n}=\bm{h} \hspace{5mm} \text{on} \hspace{3mm} \partial \Omega_{h} \hspace{7mm} \text{and} \hspace{7mm} \nabla \phi \cdot \bm{n} =0 \hspace{5mm} \text{on} \hspace{3mm} \partial \Omega.
    \centering
\end{equation}
\noindent where $\bm{n}$ denotes the outward unit vector normal to the surface \(\partial \Omega\).

\section{Finite element implementation}
\label{Sec:Numerical}
This section contains the details of the numerical implementation in a finite element setting. First, some numerical considerations are presented for the phase field in section \ref{Sec:Irreversibility}, to prevent crack healing and crack growth from tensile stresses. Afterwards, the discretisation of the problem and the formulation of residuals and stiffness matrices is addressed in section \ref{2.5}. The solution strategies considered in this paper are presented in section \ref{Sec:StaggeredSolutionScheme} and details of the Broyden-Fletcher-Goldfarb-Shanno (BFGS) algorithm for the quasi-Newton solution technique are provided in section \ref{sec:BFGS}. Finally, a scheme for selectively reducing the increment size within a monolithic approach is presented in section \ref{sec:TimeStep}. The implementation is carried out in the commercial finite element package Abaqus by means of a user element (UEL) subroutine. Abaqus2Matlab is employed to pre-process the input files \cite{AES2017}.

\subsection{Addressing irreversibility and crack growth in compression}
\label{Sec:Irreversibility}

First, a history variable field $H$ is introduced to ensure damage irreversibility,
\begin{equation} \label{growth}
    \phi_{t+\Delta t} \geq \phi_{t}
    \centering
\end{equation}

\noindent where $\phi_{t+\Delta t}$ is the phase field variable in the current time increment while $\phi_{t}$ denotes the value of the phase field on the previous increment. To ensure irreversible growth of the phase field variable, the history field must satisfy the Kuhn-Tucker conditions
\begin{equation}
    \psi_{0} - H \leq 0 \text{,} \hspace{7mm} \dot{H} \geq 0 \text{,} \hspace{7mm} \dot{H}(\psi_{0}-H)=0
    \centering
\end{equation}
\noindent for both loading and unloading scenarios. Accordingly, the history field may for a current time $t$ be written as: 
\begin{equation}
     H = \max_{\tau \in[0,t]}\psi_0( \tau). 
\end{equation}

Other approaches such as crack-sets \cite{Pham2011,Tanne2018} or penalty-based methods \cite{Wheeler2014,Gerasimov2019} have been proposed, and the treatment of the irreversibility constraint is receiving increasing attention. See the recent work by Gerasimov and De Lorenzis \cite{Gerasimov2019} for a detailed discussion and comparative studies.\\

Second, we introduce a strain energy decomposition to prevent cracking in compression. A few options have been proposed in the literature, of which the most popular ones are: the spectral tension-compression decomposition by Miehe \textit{et al.} \cite{Miehe2010a} and the volumetric-deviatoric split by Amor \textit{et al.} \cite{Amor2009}. Both are considered here but the latter will be generally adopted, unless otherwise stated. Both models have a similar intent: to maintain resistance in compression and during crack closure. In the volumetric-deviatoric split by Amor \textit{et al.} \cite{Amor2009}, the idea is that the volumetric and deviatoric strain energies can be subjected to damage but not the compressive volumetric strain energy. Thus, the strain energy can be decomposed as $\psi=g(\phi)\psi_0^{+} + \psi_0^{-}$, where
\begin{equation}
\psi_0^+=\frac{1}{2} K_n \langle tr \left( \bm{\varepsilon} \right) \rangle^2_{+} + \mu \left( \bm{\varepsilon}^{dev} : \bm{\varepsilon}^{dev} \right)
\end{equation}
\begin{equation}
\psi_0^-=\frac{1}{2} K_n \langle tr \left( \bm{\varepsilon} \right) \rangle^2_{-}
\end{equation}

\noindent where $K_n=\lambda+2\mu/n$ is bulk modulus (with $n$ being the number of dimensions of the problem), $\langle a \rangle_{\pm}=(a\pm |a|)/2$ and $\bm{\varepsilon}^{dev}=\bm{\varepsilon}-tr(\bm{\varepsilon}) \bm{I} /3$.
We follow the hybrid implementation of Ambati \textit{et al.} \cite{Ambati2015a} in considering $\psi_0^+$ in the evaluation of the history variable field $H$, therefore referring to it as $H^+$ henceforth, while considering $\psi_0$ in the displacement problem. In this regard, we emphasize that our findings relate to the incrementally linear problem resulting from the hybrid approach by Ambati \textit{et al} \cite{Ambati2015a}; the performance of monolithic quasi-Newton methods in other models remains unaddressed.

\subsection{Finite element discretisation of variational principles}
\label{2.5}

Consider (\ref{eq:weak0}), in the absence of body forces, the two-field weak form can be formulated as
\begin{equation} \label{eq:weak}
    \begin{split}
        \int_{\Omega} \left(\rho \ddot{\bm{u}} \, \delta \bm{u}+  \bm{\sigma} \delta \bm{\varepsilon} - \mathbf{b} \cdot \delta \mathbf{u} \right) \, \mathrm{d}V - \int_{\partial \Omega_{h}} \mathbf{h} \cdot \delta \mathbf{u}\, \mathrm{d}S &= 0\\[1mm]
        \int_{\Omega} \left\{ -2(1-\phi)\delta \phi \, H^+ +
        G_{c}\left(\dfrac{1}{\ell}\phi \delta \phi
        + \ell\nabla \phi \cdot \nabla \delta \phi \right) \right\}  \, \mathrm{d}V &= 0
    \end{split}
    \centering
\end{equation}

Now make use of Voigt notation and consider a plane strain solid. The displacement field \(\textbf{u}\) and the phase field \(\phi\) can be discretised as
\begin{equation}
    \mathbf{u}=\sum_{i=1}^{m} \mathbf{N}_{i}^{\mathbf{u}} \mathbf{u}_{i} \hspace{7mm} \text{and} \hspace{7mm} \phi=\sum_{i=1}^{m} N_{i} \phi_{i} 
    \centering
\end{equation}

\noindent where the shape function matrix is expressed as
\begin{equation}
    \mathbf{N}_{i}^{\mathbf{u}}=
        \begin{bmatrix}
            N_{i}   &   0  \\[0.3em]
            0   &   N_{i} 
        \end{bmatrix}
    \centering
\end{equation}
\noindent Here, \(N_{i}\) denotes the shape function associated with node \(i\), \(m\) is the total number of nodes per element, and \(\textbf{u}_{i}=\left\{u_{x}, u_{y}\right\}^{T}\) and \(\phi_{i}\) are the displacement and phase field values at node \(i\), respectively. Consequently, the corresponding derivatives can be discretised as
\begin{equation}
    \bm{\varepsilon}=\sum_{i=1}^{m} \mathbf{B}_{i}^{\mathbf{u}} \mathbf{u}_{i} \hspace{7mm} \text{and} \hspace{7mm} \nabla \phi=\sum_{i=1}^{m} \mathbf{B}_{i}^{\phi} \phi_{i}
    \centering
\end{equation}
\noindent where \(\bm{\varepsilon}=\left\{\varepsilon_{xx}, \varepsilon_{yy}, \gamma_{xy}\right\}^{T}\). Here, $\gamma$ denotes the engineering strain, such that $\gamma_{xy}=2 \varepsilon_{xy}$. Accordingly, the strain-displacement matrices associated with a given node $i$ are expressed as
\begin{equation}
    \mathbf{B}_{i}^{\mathbf{u}}=
        \begin{bmatrix}
            \partial N_i / \partial x   &   0    \\[0.3em]
            0   &    \partial N_i / \partial y    \\[0.3em]
            \partial N_i / \partial y  &   \partial N_i / \partial x 
        \end{bmatrix} \hspace{7mm} \text{and} \hspace{7mm}
        \mathbf{B}_{i}^{\phi}=
        \begin{bmatrix}
            \partial N_i / \partial x \\[0.3em]
            \partial N_i / \partial y 
        \end{bmatrix}
    \centering
\end{equation}

\noindent Considering this finite element discretisation and the fact that (\ref{eq:weak}) must hold for arbitrary values of \(\delta \textbf{u}\) and \(\delta \phi\), the discrete equation corresponding to the equilibrium condition can be expressed as the following residual with respect to the displacement field
\begin{equation}
    \mathbf{r}_{i}^\mathbf{u}=\int_{\Omega} \left[(1-\phi)^{2}+k\right] {(\mathbf{B}_{i}^{\mathbf{u}})}^{T} \bm{\sigma_{0}} \, \mathrm{d}V + \int_{\Omega} \rho {(\mathbf{N}_{i}^{\mathbf{u}})}^{T} \ddot{\bm{u}}   \, \mathrm{d}V  - \int_{\partial \Omega_{h}} {(\mathbf{N}_{i}^{\mathbf{u}})}^{T} \mathbf{h} \, \mathrm{d}S
    \centering
\end{equation}

\noindent where $k$ is a numerical parameter introduced to keep the system of equations well-conditioned. Similarly, the residual with respect to the evolution of the crack phase field can be expressed as
\begin{equation}
    r_{i}^{\phi}= \int_{\Omega} \left\{ -2(1-\phi) N_{i} \, H^+ +
    G_{c}\left[\dfrac{1}{\ell} N_{i} \, \phi
    + \ell {(\mathbf{B}_{i}^{\phi})}^{T} \nabla \phi \right] \right\} \, \mathrm{d}V
    \centering
\end{equation}

\subsection{Solution schemes}
\label{Sec:StaggeredSolutionScheme}

The Newton-Raphson method is employed to obtain the solutions for which \(\textbf{r}^{\textbf{u}}=\bm{0}\) and \(\textbf{r}^{\phi}=\bm{0}\), given the nonlinearity of the residuals. An iterative scheme is adopted to solve for the displacement $\bm{u}$ and the phase $\phi$. The tangent stiffness matrices and the mass matrix can be readily computed by taking the first derivative of the residual vectors, and read
\begin{equation}
    \begin{split}
        \mathbf{K}_{ij}^{\mathbf{u}\mathbf{u}} &= \dfrac{\partial \mathbf{r}_{i}^\mathbf{u}}{\partial \mathbf{u}_{j}} = 
        \int_{\Omega} \left[(1-\phi)^{2}+k\right] {(\mathbf{B}_{i}^{\mathbf{u}})}^{T} \mathbf{C_{0}} \, \mathbf{B}_{j}^{\mathbf{u}} \, \mathrm{d}V \\[3mm]
        \mathbf{K}_{ij}^{\phi\phi} &= \dfrac{\partial r_{i}^{\phi}}{\partial \phi_{j}} = \int_{\Omega} \left\{ \left[ 2H^+ + \dfrac{G_{c}}{\ell} \right] N_{i} N_{j} + G_{c} \ell {(\mathbf{B}_{i}^{\phi})}^{T} (\mathbf{B}_{j}^{\phi}) \right\} \, \mathrm{d}V\\[3mm]
        \mathbf{M} &= \int_{\Omega} \rho (\mathbf{N}^{\bm{u}}_i)^T\mathbf{N}^{\bm{u}}_i\, \tm{d}V.
    \end{split}
    \centering
\end{equation}

Two solution approaches are generally used to solve the phase field - displacement system: (1) solving for $\bm{u}$ and $\bm{\phi}$ simultaneously (\emph{monolithic}) or (2) solving for $\bm{u}$ and $\bm{\phi}$ separately as sequentially coupled \emph{staggered} fields. Staggered solution schemes are very robust and can overcome snap-back instabilities. However, they are not unconditionally stable and the time increment must be sufficiently small to prevent deviating from the equilibrium solution. On the other hand, monolithic implementations retain unconditional stability, enabling to use much larger time increments. Notwithstanding, the use of the more efficient monolithic schemes has been hindered by their poor performance in attaining a converged solution. We show that this Achilles' heel of monolithic solution schemes can be overcome by using quasi-Newton methods, such as the Broyden-Fletcher-Goldfarb-Shanno (BFGS) algorithm. The performance of the BFGS algorithm will be compared to that of a staggered solution scheme where convergence is assessed independently for the displacement and phase fields at the end of each increment. This widely used approach is typically referred to as one-pass or single-iteration alternating minimisation solver. The reader is referred to the recent work by Wu \textit{et al.} \cite{Wu2020a} for a comparison between the BFGS algorithm and the staggered approach employed by Bourdin \textit{et al.} \cite{Bourdin2000}, which iterates on the current phase field and displacement solutions. 

\subsection{The Broyden-Fletcher-Goldfarb-Shanno (BFGS) algorithm}\label{sec:BFGS}

Consider the following linearized system, with initial stiffness matrix $\mathbf{K}$, to be solved in an iterative manner,
\begin{equation}\label{Eq:GlobalElementSystem}
    {\begin{Bmatrix}
        \textbf{u}\\[0.3em] \bm{\phi}
    \end{Bmatrix}}_{t+\Delta t} = 
    {\begin{Bmatrix}
        \textbf{u}\\[0.3em] \bm{\phi}
    \end{Bmatrix}}_{t} -
    {\begin{bmatrix}
        \mathbf{K}^{\mathbf{u}\mathbf{u}}+\mathbf{M} & 0 \\[0.3em] 
        0 & \mathbf{K}^{\phi\phi}
    \end{bmatrix}}_{t}^{-1}
    {\begin{Bmatrix}
        \textbf{r}^{\textbf{u}}\\[0.3em] \textbf{r}^{\phi}
    \end{Bmatrix}}_{t}
    \centering
\end{equation}

In quasi-Newton methods, in contrast to standard Newton, the stiffness matrix $\mathbf{K}$ is not updated after each iteration. Instead, after a set number of iterations without convergence, an approximation of the stiffness $\tilde{\mathbf{K}}$ is introduced. This approximated stiffness matrix $\tilde{\mathbf{K}}$ satisfies the following: 
\begin{equation}
    \tilde{\mathbf{K}}\Delta\mathbf{z}= \Delta\mathbf{r}
\end{equation}
where 
\begin{equation*}
    \mathbf{z}=\begin{Bmatrix} \mathbf{u}\\[0.3em] \bm{\phi}\end{Bmatrix}
\end{equation*}
and $\Delta\mathbf{z}= \mathbf{z}_{t+\Delta t}-\mathbf{z}_t$. Likewise, $\Delta\mathbf{r}=\mathbf{r}_{t+\Delta t}-\mathbf{r}$. In the BFGS algorithm, the approximated stiffness matrix is updated in the following way: 
\begin{equation}\label{eq:BFGSfirst}
    \tilde{\mathbf{K}} = \tilde{\mathbf{K}}_t - \dfrac{(\tilde{\mathbf{K}}_t\Delta\mathbf{z})(\tilde{\mathbf{K}}_t)\Delta\mathbf{z})^T}{\Delta\mathbf{z}\tilde{\mathbf{K}}_t\Delta\mathbf{z}}+\dfrac{\Delta\mathbf{r}\Delta\mathbf{r}^T}{\Delta\mathbf{z}^T\Delta\mathbf{r}}
\end{equation}

Note that, although the non-diagonal coupling terms of the initial stiffness matrix have been dropped, see (\ref{Eq:GlobalElementSystem}), the approximation (\ref{eq:BFGSfirst}) couples the displacement and phase fields. Also, if the stiffness matrix is symmetric, the update to the approximate stiffness matrix can instead be written in terms of its inverse \cite{Matthies1979}: 
\begin{equation}
    \tilde{\mathbf{K}}^{-1} = \left(\mathbf{I}-\dfrac{\Delta\mathbf{z}\Delta\mathbf{r}^T}{\Delta\mathbf{z}^T\Delta\mathbf{r}} \right)\tilde{\mathbf{K}}_t^{-1}\left(\mathbf{I}-\dfrac{\Delta\mathbf{z}\Delta\mathbf{r}^T}{\Delta\mathbf{z}^T\Delta\mathbf{r}}\right)^{-1}+\dfrac{\Delta\mathbf{z}\Delta\mathbf{z}^T}{\Delta\mathbf{z}^T\Delta\mathbf{r}},
\end{equation}
which offers significant computational savings and retains symmetry and positive definiteness, if such was already present. The BFGS algorithm has been implemented in most commercial finite element packages (such as Abaqus), often in conjunction with a line search algorithm.

\subsection{Convergence criteria}
\label{sec:Convergence}

The standard convergence criteria available in Abaqus are used for both the monolithic and staggered solution schemes without any modification. Hence, both a residual control and a solution correction control have to be met to achieve convergence. Regarding the former, the largest residual in the balance equations $r^\alpha_{max}$ must be equal or smaller than the product of a tolerance constant $R^\alpha_n$ with an overall time-averaged flux norm for the solution $\tilde{q}^\alpha$:
\begin{equation}\label{eq:ResidualTol}
   r^\alpha_{max} \leq  R^\alpha_n \tilde{q}^\alpha
\end{equation}

We do not deviate from Abaqus default recommendations and consider a magnitude for the tolerance of $R^\alpha_n=0.005$. If the inequality (\ref{eq:ResidualTol}) is satisfied, convergence is accepted if the largest correction to the solution, $c^\alpha_{max}$, is also small compared to the largest incremental change in the corresponding solution variable, $\Delta a_{max}^\alpha$,
\begin{equation}
c^\alpha_{max} \leq C^\alpha_n \Delta a_{max}^\alpha     
\end{equation}

Here, $C^\alpha_n$ denotes the magnitude of the convergence tolerance; as by default in Abaqus, we consider $C^\alpha_n=0.01$. The residual-based and solution-based criteria are equally employed for the displacement field $(\alpha=u)$ and the phase field $(\alpha=\phi)$.

\subsection{Incrementation control scheme}\label{sec:TimeStep}

As it will be shown below, phase field fracture problems can be solved employing very large load increments when combining quasi-Newton algorithms with monolithic solution schemes. However, a sudden change in the material response (such as a large force drop due to unstable fracture) may be best captured by using small time increments. To benefit from the use of large time increments while resolving sudden changes in material behaviour we suggest the use of the following adaptive step scheme:
\begin{center}
    \begin{tabular}{l}\hline
        For any integration point $i$:\\
         \hspace{0.8cm} If $\phi_i<0.7$ \& $\Delta\phi_i \geq 0.5$\\
         \hspace{1.4cm} Re-start load increment reducing its size by 90\%: $\Delta t_1=0.1 \Delta t_0$. \\
         \hspace{0.8cm} End if\\\hline
    \end{tabular}
\end{center}

Accordingly, smaller increments are used when there is a significant increase in damage in a material point that was not already highly damaged. This allows us to generalise the present time stepping scheme to case studies where the crack is initially introduced by prescribing $\phi=1$. Obviously, the above criterion must be restrained such that it does not happen continuously. For simplicity, it has been restricted to happen only once during a given simulation. This criterion is particularly useful for unstable cracking problems, where large increments can be adopted until the onset of cracking and the load increment decreases to adequately resolve the fracture event.

\section{Results}
\label{Sec:Results}

We proceed to showcase, via numerical experiments, the potential of the method in attaining convergence and reducing computation times in a wide range of problems. First, quasi-static fracture is considered in Section \ref{sec:QuasiStatic}, including both stable and unstable cracking. Secondly, in Section \ref{sec:Fatigue}, we show that the quasi-Newton method can enable cycle-by-cycle phase field fatigue calculations that are computationally prohibitive for staggered schemes. Finally, the capabilities of the method are also demonstrated for the case of dynamic crack branching in Section \ref{sec:Dynamics}.

\subsection{Quasi-static fracture}\label{sec:QuasiStatic}

Two paradigmatic benchmarks will be addressed under quasi-static loading conditions, the fracture of a cracked square plate under: (1) uniaxial tension, and (2) shear.

\subsubsection{Cracked square plate subjected to uniaxial tension}

\begin{figure}[H]
    \begin{subfigure}[h]{0.32\textwidth}
    \includegraphics[width=0.9\linewidth]{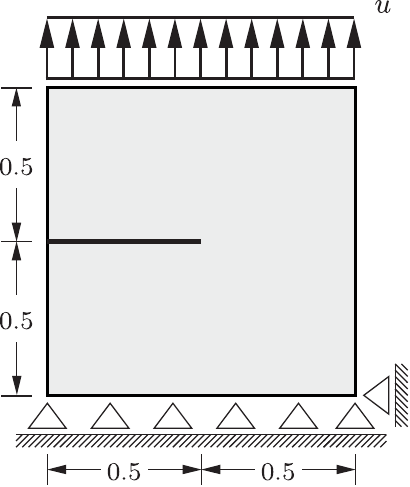}
    \caption{ }
    \label{fig:SENTGeometry}
    \end{subfigure}
    \begin{subfigure}[h]{0.32\textwidth}
    \includegraphics[width=0.95\textwidth]{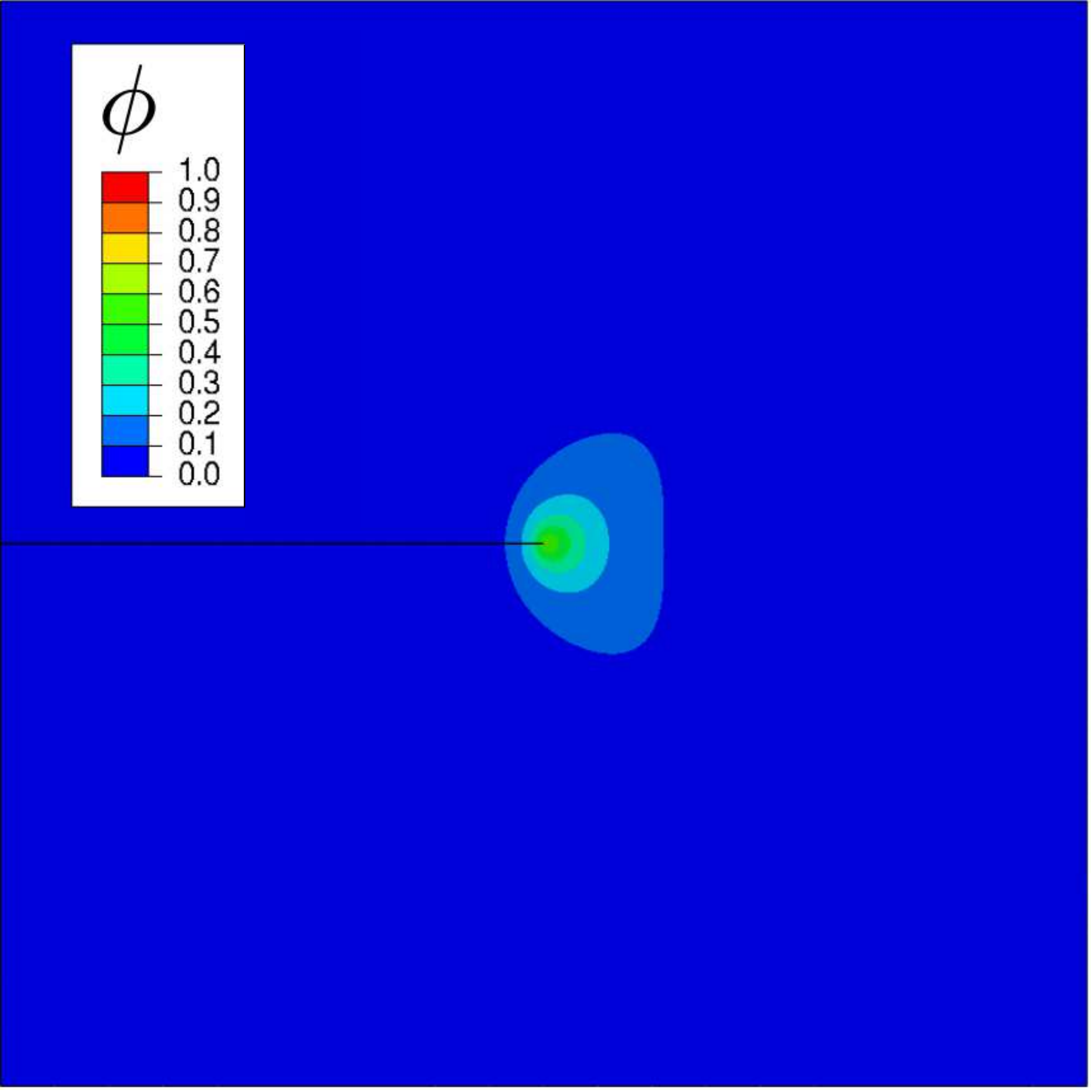}
    \caption{}
    \label{fig:SENTCrack1}
    \end{subfigure}
    \begin{subfigure}[h]{0.32\textwidth}
    \includegraphics[width=0.95\textwidth]{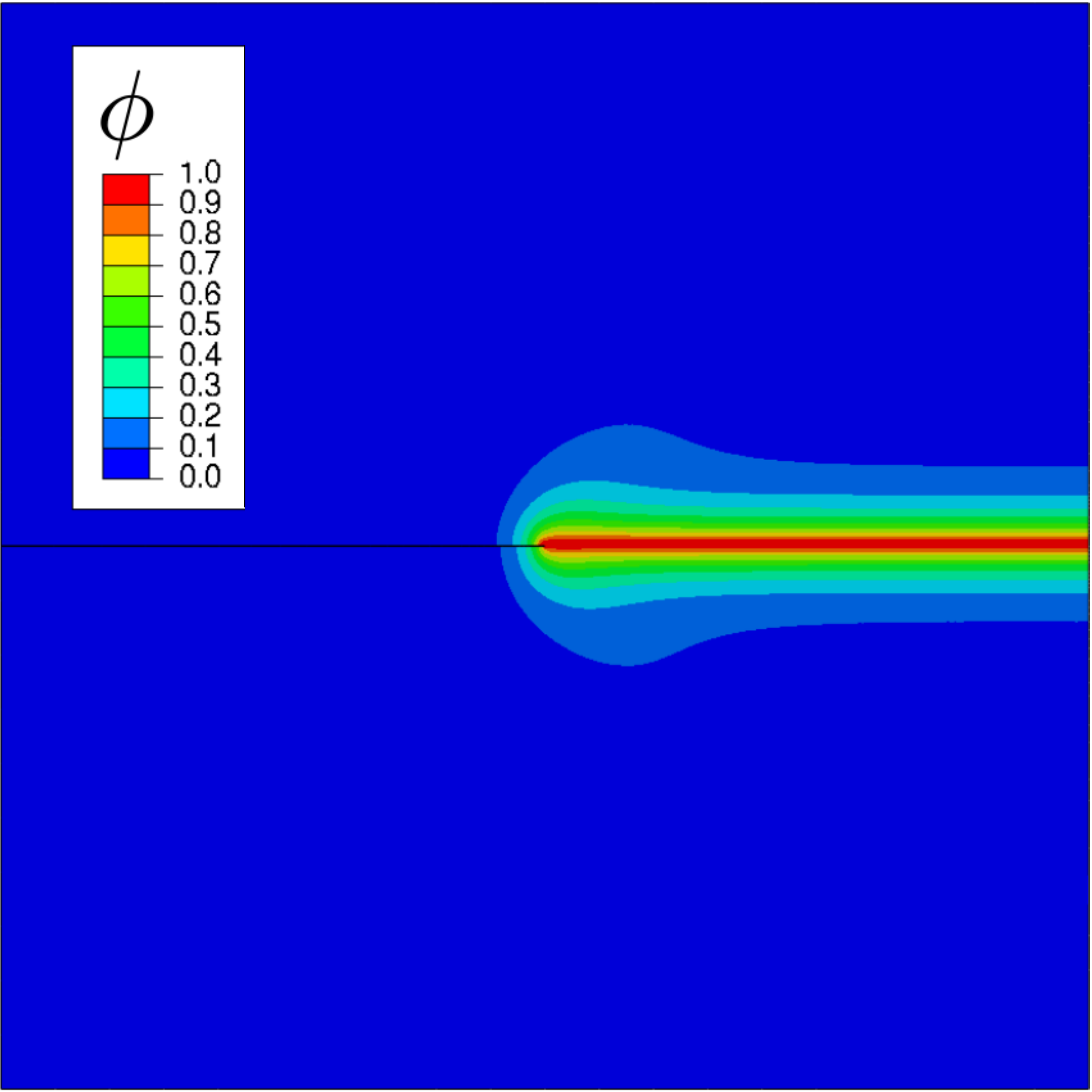}
    \caption{ }
    \label{fig:SENTCrack2}
    \end{subfigure}
    \caption{Single-edged notched tension specimen: (a) Dimensions (in mm) and loading configuration. Crack patterns for a remote displacement of (b) $u=5.90\times 10^{-3}$mm and (c) $u=5.93\times 10^{-3}$mm. The complete extent of crack growth is captured within a single increment.}
    \label{fig:SENTintro}
\end{figure}

We consider first the case of unstable crack growth in a linear elastic specimen under monotonic loading, as exemplified by the mode I fracture of the single-edge notched tension (SENT) specimen sketched in Fig. \ref{fig:SENTGeometry}. This paradigmatic example has been widely used since the early works by Miehe and co-workers \cite{Miehe2010a}. Loading conditions and specimen dimensions (in mm) are shown in Fig. \ref{fig:SENTGeometry}. Material properties are $E=210$ GPa, $\nu = 0.3$, $\ell=0.024$ mm, and $G_c=2.7$ J/$\tm{mm}^2$. We discretise the model with linear quadrilateral elements, with the characteristic element size along the extended crack path, $h_e$, being at least 6 times smaller than the phase field length scale $\ell$. Cracking is unstable, with damage extending through the crack ligament instantaneously. With the quasi-Newton monolithic scheme, the unstable growth is captured within a single increment without convergence problems. This is illustrated by means of phase field contours in Figs. \ref{fig:SENTCrack1} and \ref{fig:SENTCrack2}.

\begin{figure}[H]
    \centering
    \begin{subfigure}[h]{0.9\linewidth}
    \centering
    \includegraphics[width=0.75\textwidth]{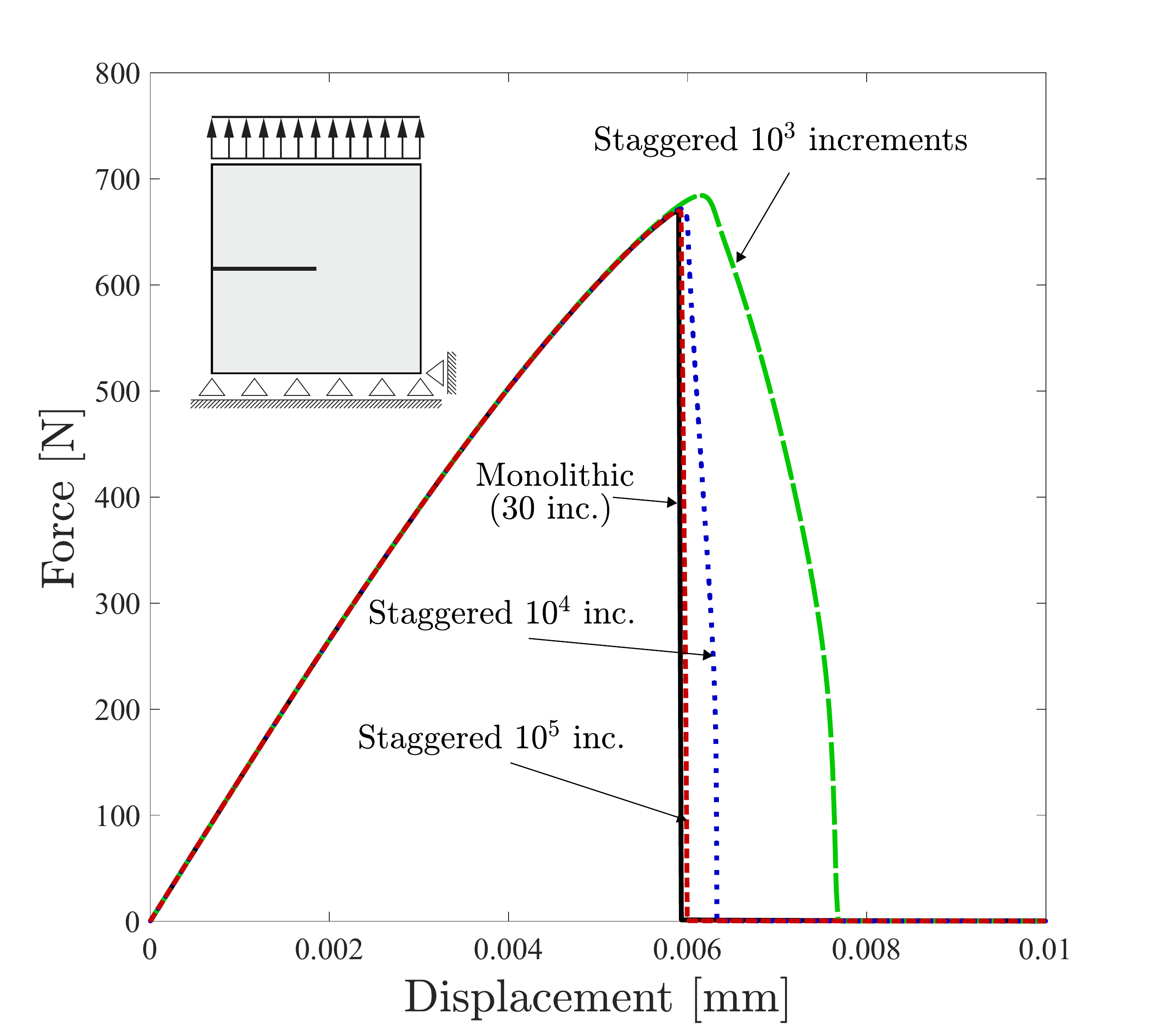}
    \caption{ }
    \label{fig:SENTFD}    
    \end{subfigure}
       \begin{subfigure}[h]{0.9\linewidth}
       \centering
    \includegraphics[width=0.75\textwidth]{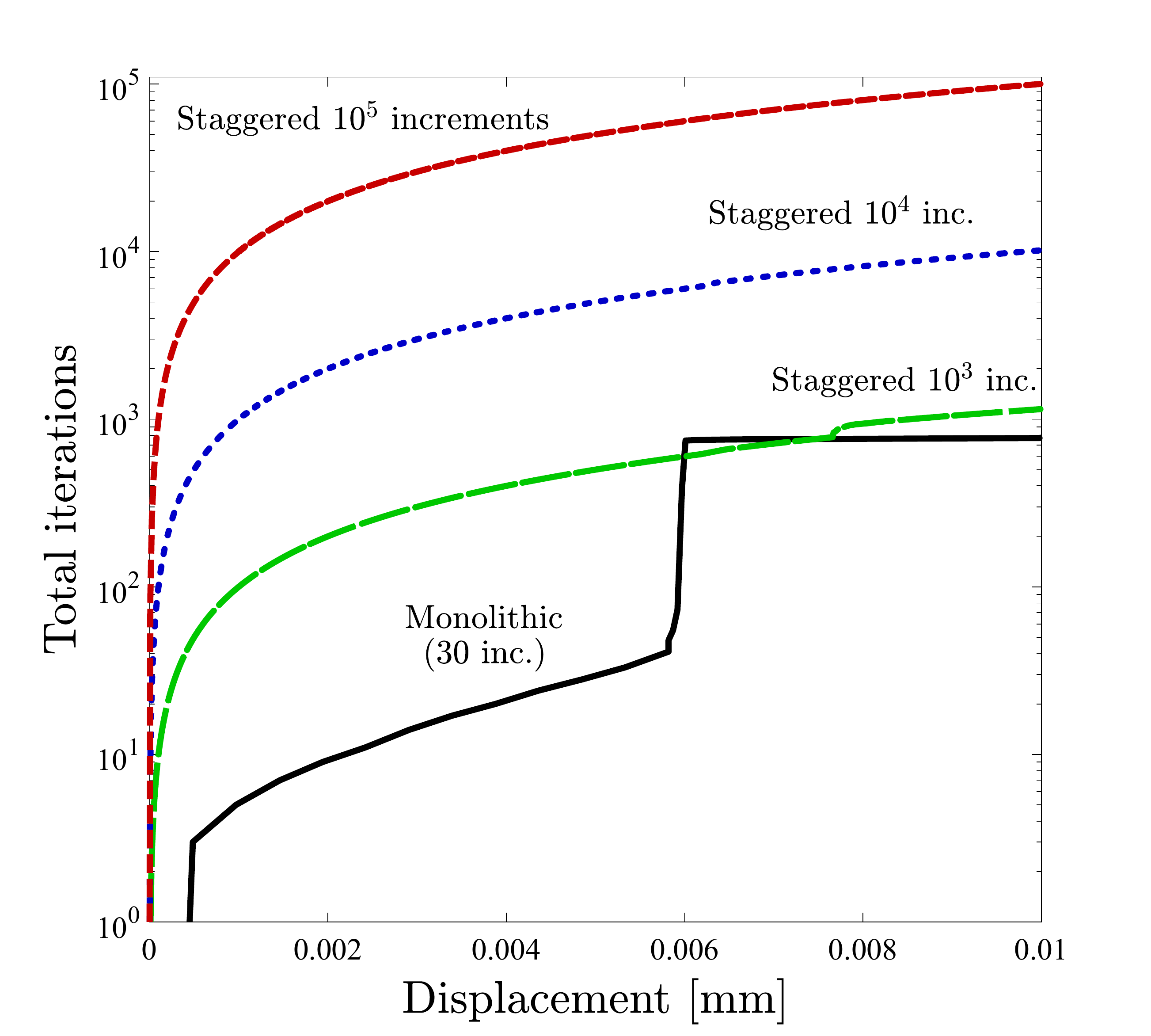}
    \caption{}
    \label{fig:SENTCummul}    
    \end{subfigure} 
\caption{Single-edged notched tension specimen: (a) force versus displacement curves, and (b) cumulative number of iterations.}
\end{figure}

The force versus displacement curve obtained is shown in Fig. \ref{fig:SENTFD}. The results computed with a staggered solution scheme are also shown for selected values of the total number of load increments employed. It can be clearly seen that the staggered solution is sensitive to the increment size, and recovering the monolithic solution requires at least $10^5$ increments. This is in clear contrast with the 30 increments employed to obtain the monolithic result. Both the staggered and monolithic computations require a large number of iterations to achieve convergence during the critical increments. The cumulative number of iterations is shown as a function of the applied displacement in Fig. \ref{fig:SENTCummul}. Results reveal that reproducing the accurate monolithic result within a staggered scheme requires using a number of iterations that is two orders of magnitude larger. The quasi-Newton monolithic calculation is faster than the coarsest staggered time-stepping, which leads to $\sim$20\% errors in the computation of the critical displacement. 

Computation times for different discretisations are reported in Table \ref{tab:SENTTimes}. Selected mesh densities are considered, as illustrated by the size of the characteristic element length along the crack path, $h_e$ and the number of degrees-of-freedom (DOFs). It can be seen that the trends persist across mesh densities; the quasi-Newton monolithic implementation presented is roughly 100 times faster than the widely used staggered scheme. These massive differences in computation times are mainly due the reduced number of increments employed in the monolithic case. 

\begin{table}[H]
    \centering
    \begin{tabular}{ c| c  c c c}
   &\multicolumn{4}{c}{CPU hours}\\ \hline
    \multirow{2}{*}{Mesh size} & $\ell/h_e=6$ & $h_e/\ell=9$ & $h_e/\ell=12$ &  $h_e/\ell=18$ \\ 
   &25908 DOFs & 47376 DOFs & 74697 DOFs & 152772 DOFs \\\hline
    Monolithic & 0.31 & 0.80 & 1.79 & 3.41 \\
    Staggered & 31.6 & 60.17 & 87.47 & 187.90 \\\hline 
 
    \end{tabular}
    \caption{Single-edged notched tension specimen. Computation times as a function of the mesh size. The staggered computations correspond to the $10^5$ increments case, which is one that exhibits a comparable accuracy to the monolithic result.}
    \label{tab:SENTTimes}
\end{table}

\begin{figure}[H]
    \centering
    \includegraphics[width=\textwidth]{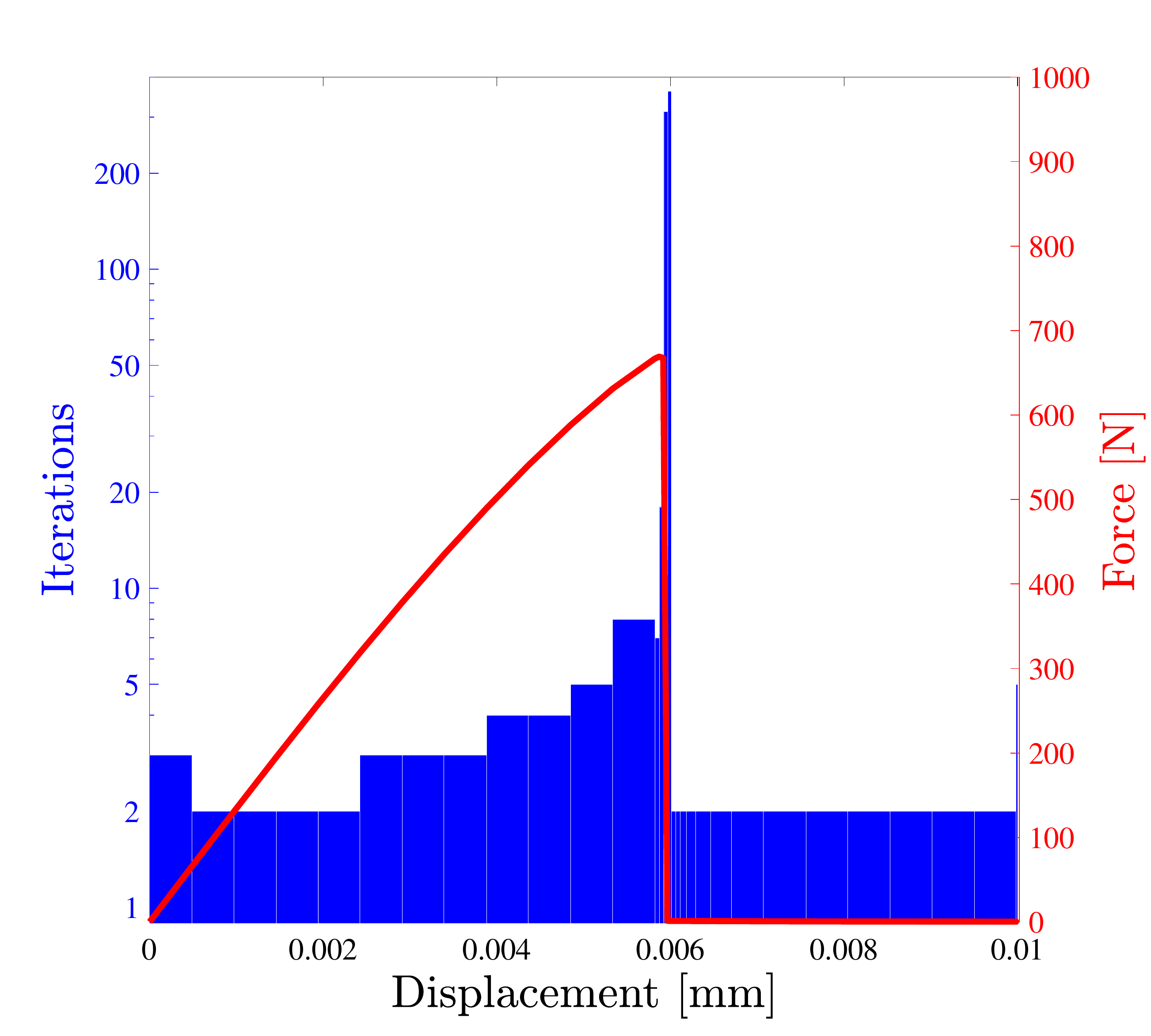}
    \caption{Single-edged notched tension specimen. Number of iterations per increment, with the force versus displacement curve superimposed. When the critical point in the simulation is reached, the increment size is drastically reduced, such that the response of the system is accurately captured with only 30 load increments.} 
    \label{fig:SENTBars}
\end{figure}

The adaptive time stepping scheme presented in Section \ref{sec:TimeStep} allows to accurately capture the unstable response while resolving with large time increments outside of the cracking time frame. This is illustrated in Fig. \ref{fig:SENTBars}, where the increment size and their required iterations are given, along with the development of the force-displacement curve. Large load increments are initially used, which require only a few iterations to converge. When cracking takes place, the algorithm drastically reduces the increment size to accurately capture unstable crack growth. As by default in Abaqus, the increment size increases when few iterations are needed to achieve convergence, enabling to recover large loading steps towards the end of the computation. 

In the context of monolithic quasi-Newton, the same force versus displacement response can be achieved without the new adaptive time stepping scheme by using a reference load increment that is approximately 85\% smaller than the one employed in Fig. \ref{fig:SENTBars}. I.e., even without adaptive time stepping, quasi-Newton calculations are 20-40 times faster than staggered ones.

\subsubsection{Single edge notched shear test}

\begin{figure}[H]
    \begin{subfigure}[h]{0.45\textwidth}
      \begin{flushleft}
    \flushleft
    \includegraphics[width=\textwidth]{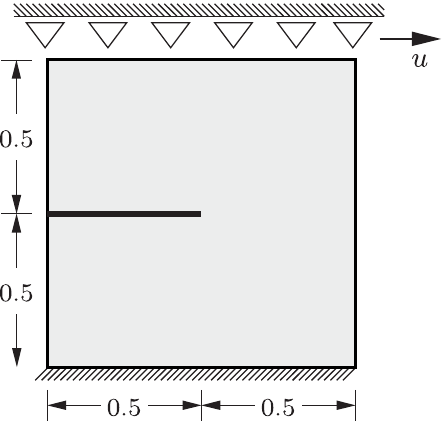}
     \label{fig:S-Geometry}
    \caption{ }
     \end{flushleft}    
     \end{subfigure}
        \begin{subfigure}[h]{0.55\textwidth}
     \begin{flushright}
     \includegraphics[width=0.8\textwidth]{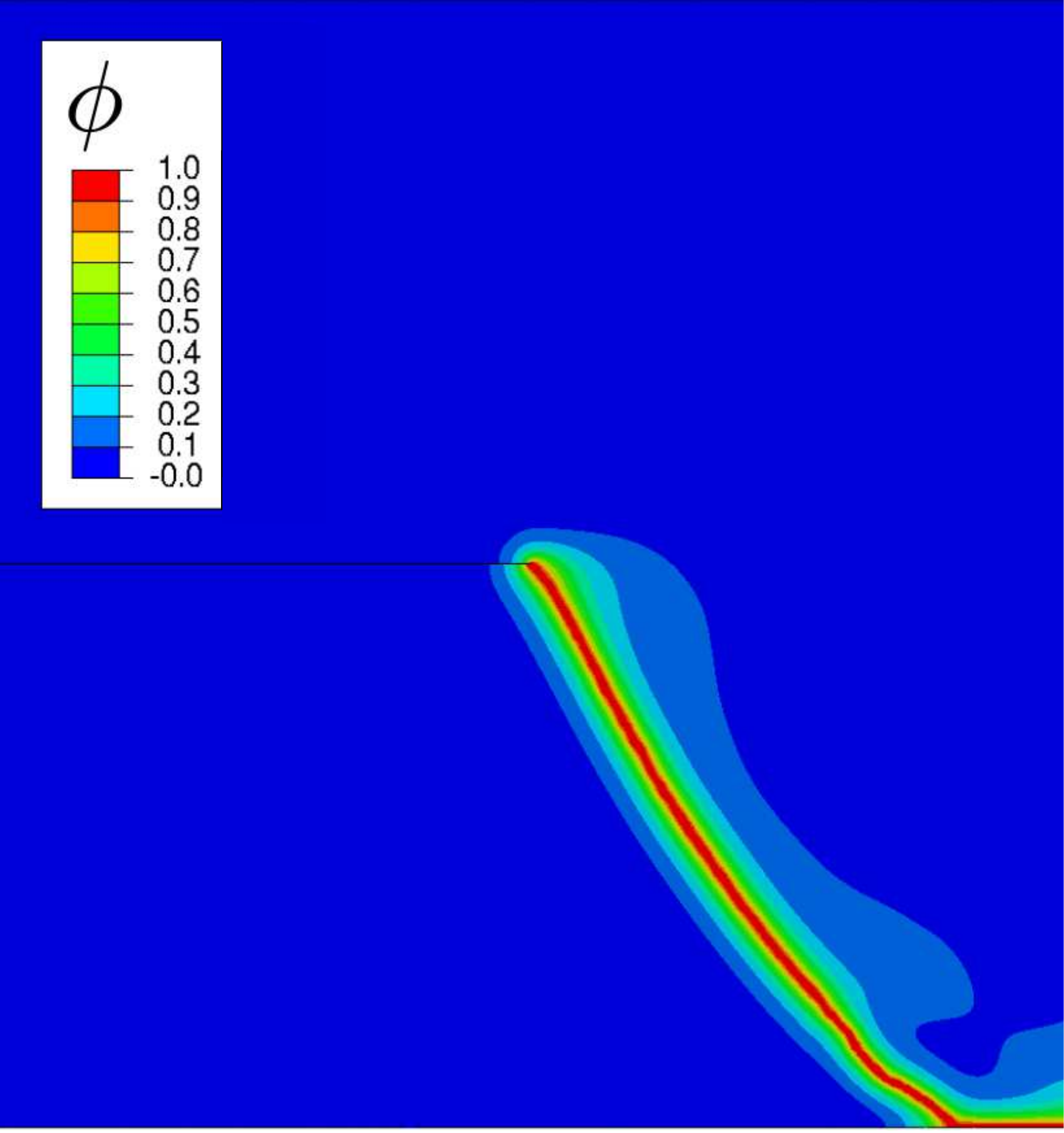}
    \label{fig:ShearCrack}
    \caption{ }
        \end{flushright} 
     \end{subfigure}
    \caption{Single-edged notched shear specimen: (a) Dimensions (in mm) and loading configuration, and (b) crack trajectory.}
    \label{fig:ShearIntro}
\end{figure}

The performance of the monolithic quasi-Newton scheme presented is now assessed in the context of stable crack growth. The same specimen dimensions and material properties as in the previous case study are employed but the sample is now subjected to remote shear loading; see Fig.  \ref{fig:ShearIntro}a. The mixed-mode crack tip conditions lead to crack deflection towards the lower part of the sample. The resulting crack trajectory, shown in Fig. \ref{fig:ShearIntro}b, agrees with results shown in the literature using the volumetric-deviatoric split - see, e.g. Ref. \cite{Ambati2015}.

\begin{figure}[H]
    \centering
    \begin{subfigure}[h]{0.9\linewidth}
    \includegraphics[width=0.75\textwidth]{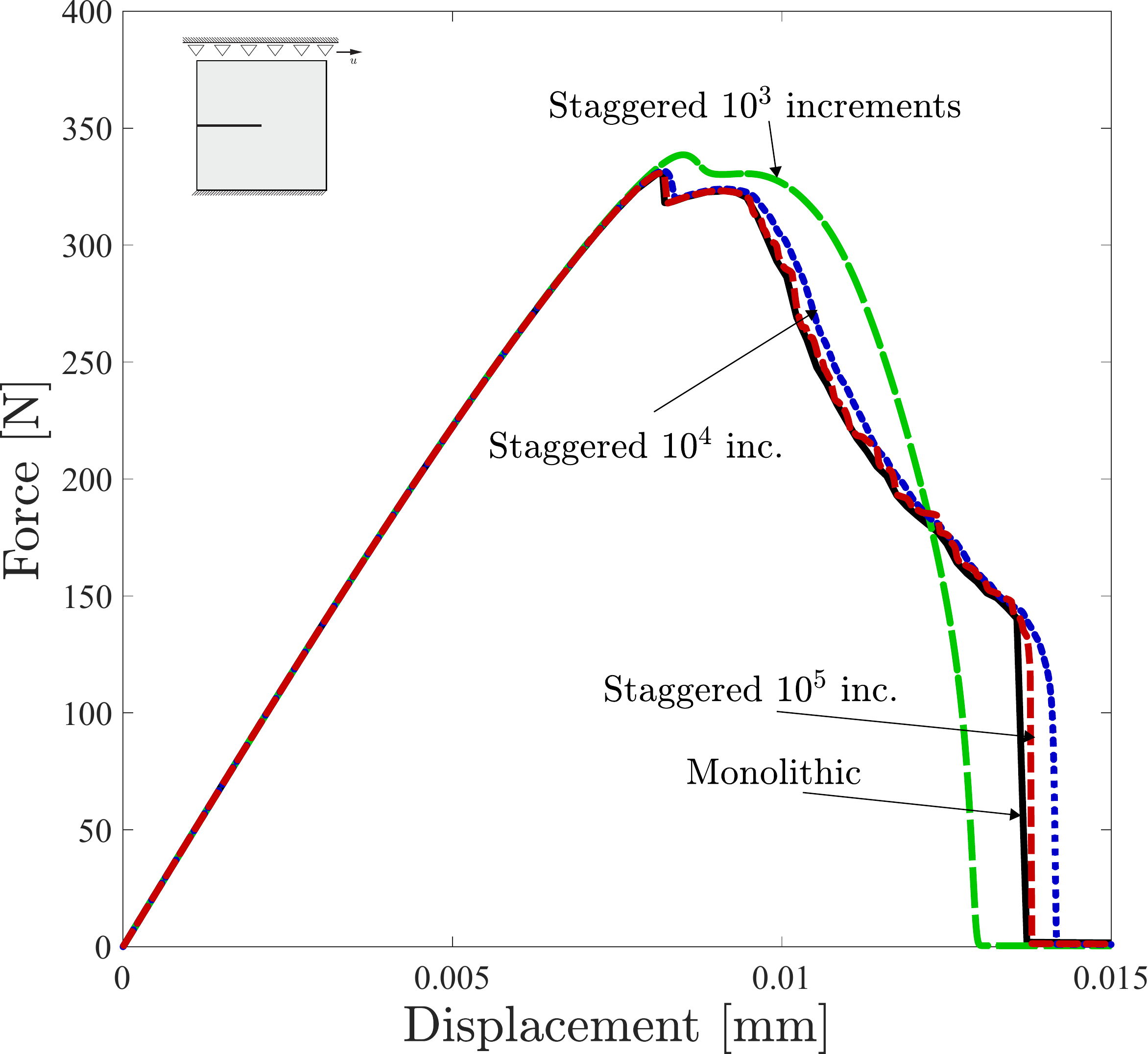}
    \caption{ }
    \label{fig:SENSFD}    
    \end{subfigure}
       \begin{subfigure}[h]{0.9\linewidth}
    \includegraphics[width=0.75\textwidth]{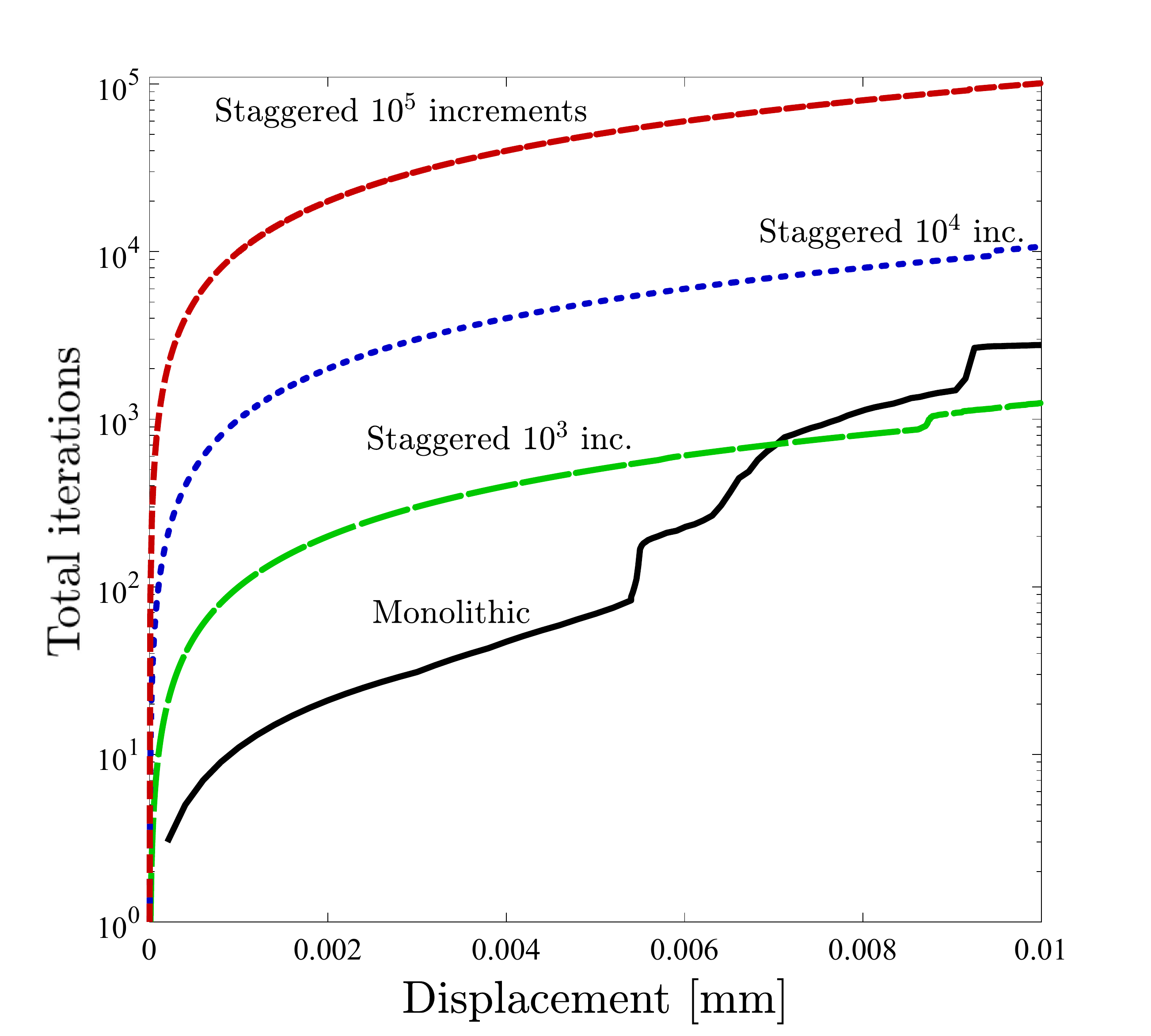}
    \caption{}
    \label{fig:SENSCummul}    
    \end{subfigure} 
\caption{Single-edged notched shear specimen: (a) force versus displacement curves, and (b) cumulative number of iterations.}
\end{figure}

The force versus displacement curves obtained with both monolithic and staggered schemes are shown in Fig. \ref{fig:SENSFD}. As with the tension case, reproducing the accurate monolithic result with the staggered implementation requires using a very large number of increments, $10^5$ or more. When using a smaller number, such as with the case of $10^3$ increments, the force versus displacement result deviates substantially from the monolithic one and a different crack path is predicted. Differences in the total number of iterations between the staggered and monolithic results are smaller than in the unstable mode I crack growth example but remain very significant. As shown in Fig. \ref{fig:SENSCummul}, the total number of iterations needed to obtain an accurate result with the staggered approach is roughly two orders of magnitude larger than in the quasi-Newton monolithic analysis.

Computation times for an increasingly refined mesh are compiled in table \ref{tab:SENSTimes}. As in the previous example, only the accurate staggered solution with $10^5$ increments is considered. In general, the computational cost is 10 times smaller in the monolithic quasi-Newton case.

\begin{table}[H]
    \centering
    \begin{tabular}{ c| c  c c c}
   &\multicolumn{4}{c}{CPU hours}\\ \hline
    \multirow{2}{*}{Mesh size} & $\ell/h_e=6$ & $h_e/\ell=9$ & $h_e/\ell=12$ &  $h_e/\ell=18$ \\ 
    & 58518 DOFs & 128451 DOFs & 222111 DOFs& 386112 DOFs \\\hline
    Monolithic & 2.02 & 6.56 & 11.62 & 46.60 \\
    Staggered & 74.85 & 159.50 & 272.25 & 469.48 \\\hline
    \end{tabular}
    \caption{Single-edged notched shear specimen. Computation times as a function of the mesh size. The staggered computations correspond to the $10^5$ increments case, which is one that exhibits a comparable accuracy to the monolithic result.}
    \label{tab:SENSTimes}
\end{table}

\begin{figure}[H]
    \centering
    \includegraphics[width=\textwidth]{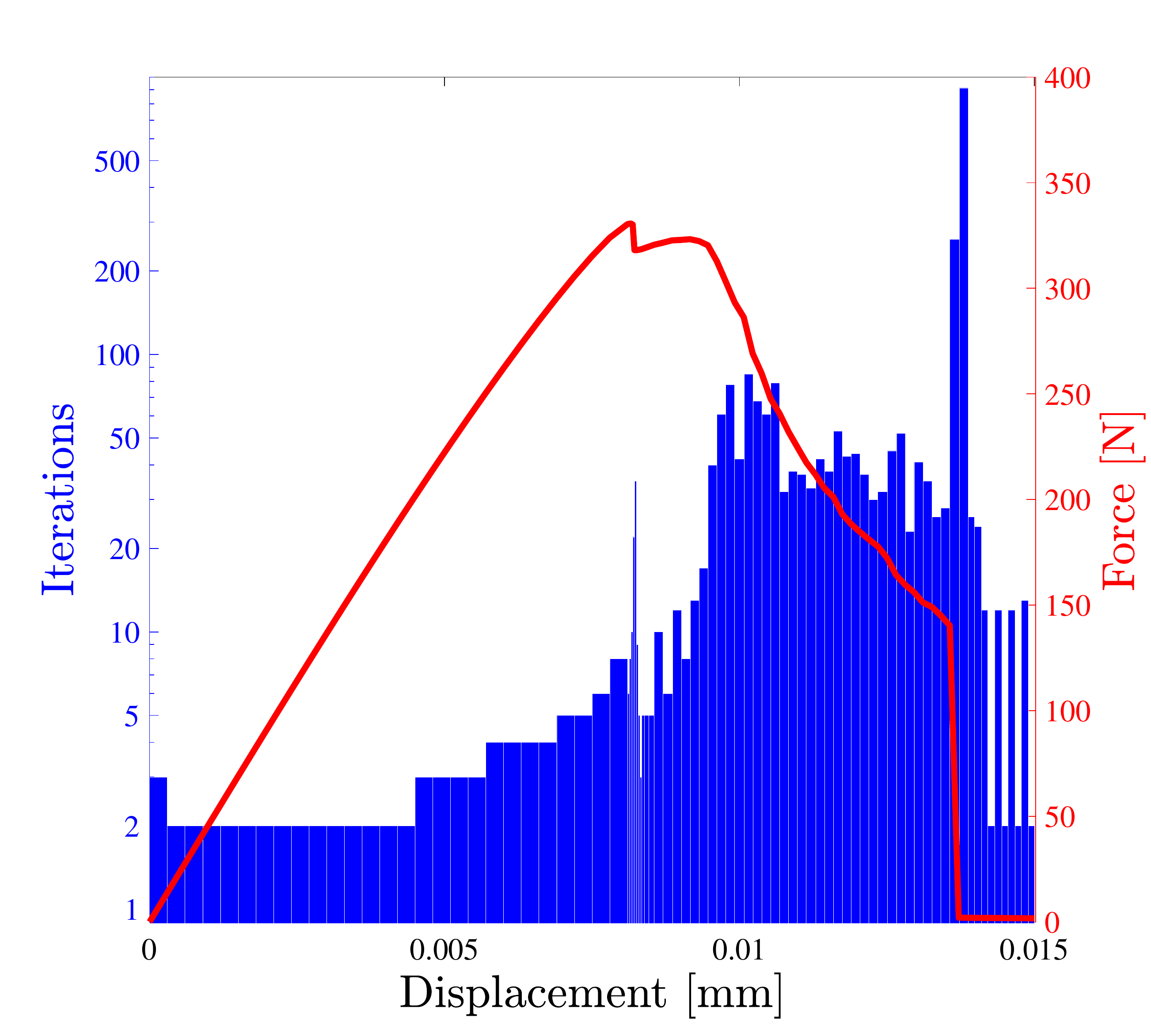}
    \caption{Single-edged notched shear specimen. Number of iterations per increment, with the force versus displacement curve superimposed.} 
    \label{fig:SENSBars}
\end{figure}

Fig. \ref{fig:SENSBars} shows the performance of the present monolithic quasi-Newton scheme, with the adaptive time stepping defined in Section \ref{sec:TimeStep}. The bar plot shows the number of iterations as a function of the remote displacement, together with the force versus displacement response. The performance is not as impressive as for the mode I unstable crack growth example but it still leads to substantial computational gains and remains useful for providing a well-timed transition from large to small increments.

Lastly, it should be noted that the method is very robust. Convergence is attained in all cases without the need for any viscous dissipation parameters. In both the tension and shear boundary value problems the monolithic implementation based on the conventional Newton method fails to converge, even for $10^5$ increments. This also holds true when combining the conventional Newton method with a line search algorithm.

\subsection{Phase field fatigue}
\label{sec:Fatigue}

We proceed to investigate the effectiveness of the monolithic quasi-Newton solution approach within the emerging field of phase field fatigue modelling. We base our investigation on the framework that Carrara \textit{et al.} \cite{Carrara2020} have very recently presented. First, a brief overview of the fatigue model is presented. Our implementation is then validated with the results by Carrara and co-workers \cite{Carrara2020}. Finally, we show that the present quasi-Newton monolithic implementation drastically outperforms staggered approaches, which are too computationally expensive for cycle-by-cycle fatigue analyses.

\subsubsection{Theoretical framework}

Consider the framework presented in Section \ref{Sec:Theory} under quasi-static conditions. A fatigue degradation function $f(\overline{\alpha}(t))$ can be introduced, which depends upon a cumulative history variable $\overline{\alpha}$ \cite{Alessi2018c,Carrara2020}. Accordingly, the variation of the internal work reads: 
\begin{equation} \label{eq:fatigueIntWork}
    \delta W_{int} =  \int_{\Omega} \left\{ \bm{\sigma} \delta \bm{\varepsilon} -2(1-\phi)\delta \phi \, \psi_{0}(\bm{\varepsilon}) +
    f(\overline{\alpha}(t))G_{c}\left(\dfrac{1}{\ell}\phi \delta \phi
    + \ell\nabla \phi \cdot \nabla \delta \phi \right) \right\}  \, \mathrm{d}V
    \centering
\end{equation}
\noindent The choices of $f$ and $\overline{\alpha}$ are of utmost importance in capturing the physics of fatigue damage. Since our aim is to investigate the performance of a new solution methodology, we restrict attention to one of the simplest choices proposed by Carrara \textit{et al.} \cite{Carrara2020}. The cumulative history variable $\overline{\alpha}$ is assumed to be independent of the mean load and takes the form: 
\begin{equation}\label{eq:Alpha}
    \overline{\alpha}(t)= \int_0^t \theta(\alpha \dot{\alpha})|\dot{\alpha}|\,\tm{d}\tau,
\end{equation}
\noindent where $\tau$ is the pseudo-time and $\theta(\alpha \dot{\alpha})$ is the Heaviside function. Thus, $\overline{\alpha}$ only grows during loading. The fatigue history variable $\alpha$ must represent the loading condition in the solid. For simplicity, the choice of $\alpha=g(\phi)\psi_0$ is made. Finally, the fatigue degradation function characterises the sensitivity of the fracture energy to the number of cycles. Here, we adopt a function that vanishes asymptotically: 
\begin{equation}
    f(\overline{\alpha}(t)) = \begin{cases} \hspace{1cm}1 &  \tm{if }\hspace{0.25cm} \overline{\alpha}(t) \leq \alpha_T \\[0.4cm] \left(\dfrac{2\alpha_T}{\overline{\alpha}(t)+\alpha_T}\right) &  \tm{if }\hspace{0.25cm} \overline{\alpha}(t) \leq \alpha_T \end{cases}.
\end{equation}
Here, $\alpha_T$ represents a threshold value, below which the fracture energy remains unaffected. The extension of the finite element implementation described in Section \ref{Sec:Numerical} to incorporate $f(\overline{\alpha}(t))$ is straightforward and will not be detailed here.

\subsubsection{Verification}

We mimic the first benchmark study by Carrara \textit{et al.} \cite{Carrara2020}. A single-edge notched tension specimen like the one described in Fig. \ref{fig:SENTGeometry} is subjected to cyclic loading with a load ratio of $R=-1$ (equal compression and tension loads). As in the original study, material parameters are chosen as $E=210$ GPa, $\nu=0.3$, $Gc = 2.7$ N/mm, $\alpha_T=56.25\, \tm{N/mm}^2$ and $\ell=0.004$ mm. The loading amplitude is of $0.002$ mm and the characteristic element size in the fracture zone is $h_e/\ell=5$. The crack extension $a$ is computed as a function of the number of cycles $N$ for three different decompositions of the strain energy density. Namely, the standard isotropic one, the volumetric-deviatoric split \cite{Amor2009}, and the spectral tension-compression decomposition \cite{Miehe2010}. The results obtained in Ref. \cite{Carrara2020} are shown superimposed using symbols. It should be noted that Carrara and co-workers \cite{Carrara2020} employ a staggered solution scheme that iterates until convergence to the monolithic solution using the energy-based convergence criterion presented in \cite{Ambati2015}.

\begin{figure}[H]
\centering
\includegraphics[width=\textwidth]{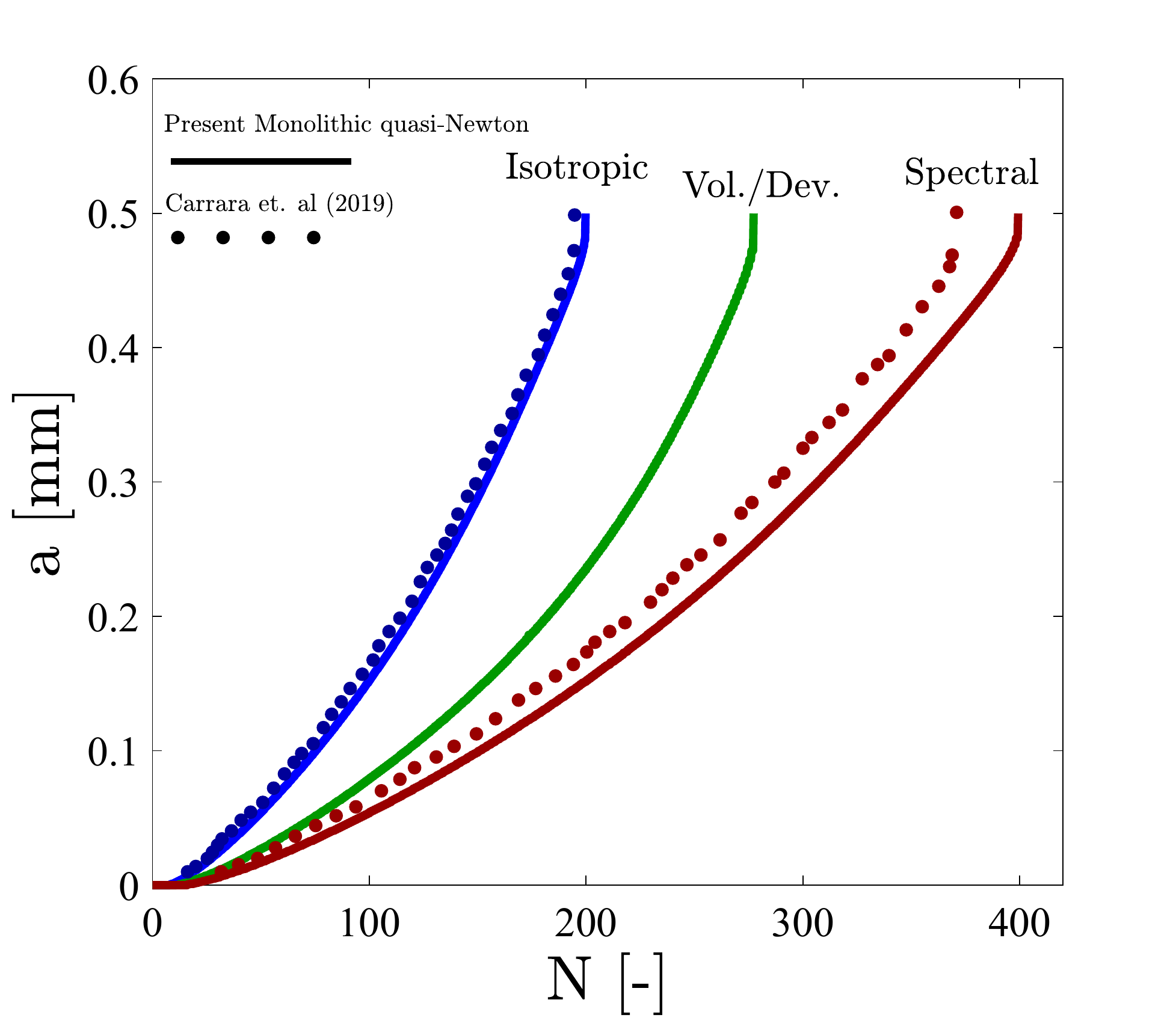}
\label{fig:StrainSplits}
\caption{Fatigue cracking of a single-edge notched tension specimen. Crack length $a$ versus number of cycles $N$ for different strain energy density splits. Solid lines are the results obtained with our quasi-Newton solution while symbols denote the results reported in Ref. \cite{Carrara2020}.}
\end{figure}

A relatively good agreement is observed for the predictions with the isotropic and spectral decomposition models. As discussed by Carrara \textit{et al.} \cite{Carrara2020}, the the spectral split requires roughly twice as many cycles to trigger complete fracture, as compressive loading cycles do not contribute to damage. Interestingly, we note that our quasi-Newton monolithic implementation of the volumetric-deviatoric split convergences without problems until final fracture.

\subsubsection{Performance of quasi-Newton in phase field fatigue}

\begin{figure}[H]
    \centering
    \includegraphics[width=\textwidth]{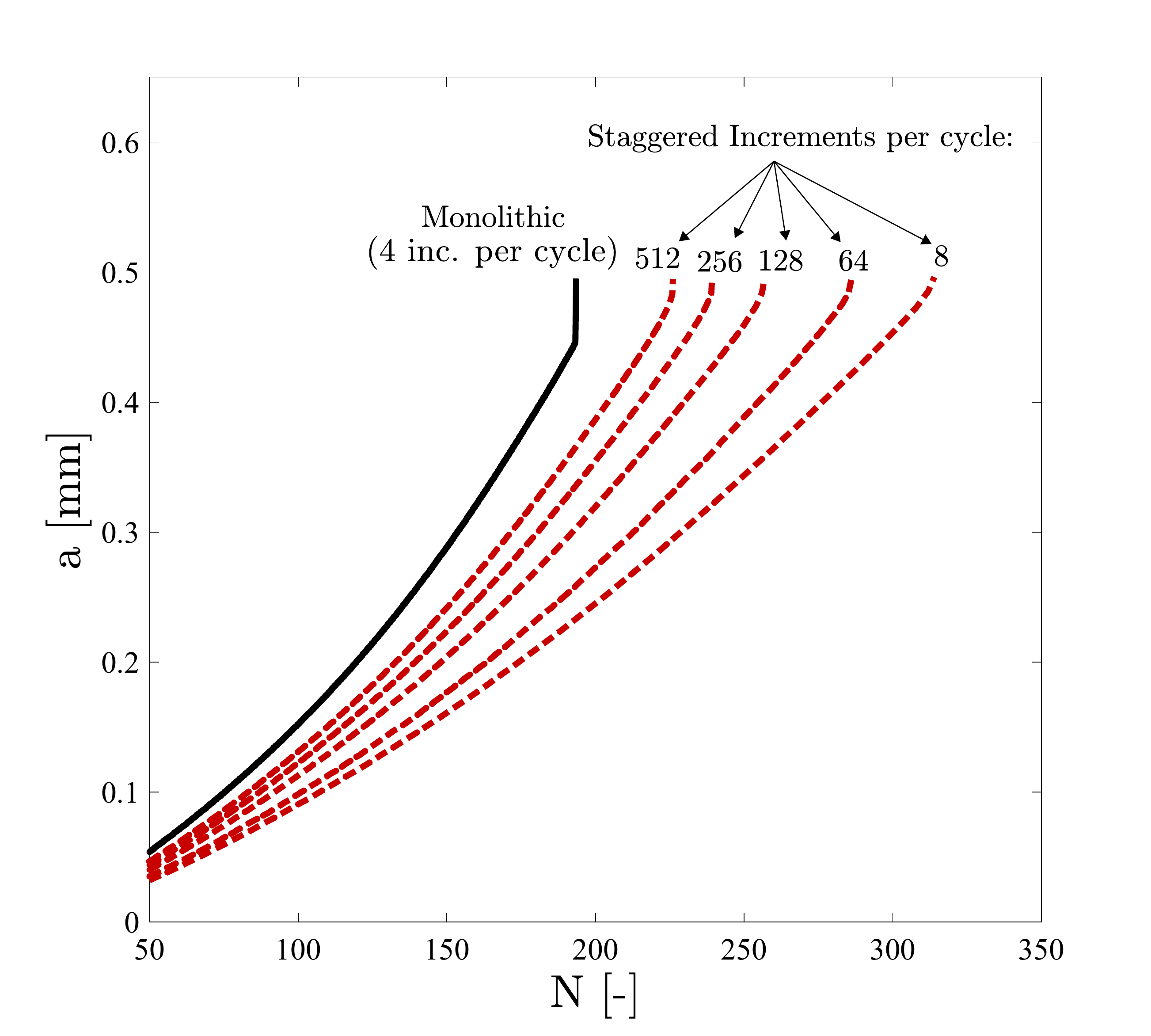}
    \caption{Fatigue cracking of a single-edge notched tension specimen. Crack length $a$ versus number of cycles $N$ for quasi-Newton monolithic and staggered implementations. A very large number of increments is needed in the staggered scheme to reproduce the quasi-Newton monolithic prediction.}    
    \label{fig:FatigueStagger}
\end{figure}

We proceed to evaluate the performance of the quasi-Newton approach. As in the other cases, we compare against the single iteration staggered approach presented in section \ref{Sec:StaggeredSolutionScheme}. For the sake of simplicity, we restrict our attention to the isotropic case, not considering any decomposition of the strain energy density. The results obtained with both staggered and monolithic approaches are shown in Fig. \ref{fig:FatigueStagger}. It can be clearly seen that a very large number of increments per cycle is needed in the staggered case, as convergence towards the monolithic solution is slow. 

The very large number of increments needed has an immediate impact on the computation times, as listed in Table \ref{tab:FatigueStagger}. The most precise staggered calculation, which is far from the equilibrium result, requires computation times that are 5 times larger than the monolithic case. From Fig. \ref{fig:FatigueStagger}, it seems likely that more than 1000 increments per cycle will be needed to approximate the monolithic result. This implies that mid, high and very high cycle fatigue problems cannot be addressed with staggered schemes; the use of quasi-Newton methods could open new horizons in phase field fatigue analyses. 

\begin{table}[H]
    \centering
    \begin{tabular}{ r| c | c  c  c  c c }
    Solutions strategy & Monolithic & \multicolumn{5}{c}{Staggered}\\\hline
    Increments per cycle     & 4& 8 & 32 & 64 & 128& 256   \\ \hline
      CPU hours   &  14.85& 3.24 & 16.52 & 20.29 & 34.30 & 73.98 \\\hline
    \end{tabular}
    \caption{Fatigue cracking of a single-edge notched tension specimen. Computation times for quasi-Newton monolithic and staggered approaches.}
    \label{tab:FatigueStagger}
\end{table}

\subsection{Dynamic results}
\label{sec:Dynamics}

\begin{figure}[H]
    \centering
    \includegraphics[width=\textwidth]{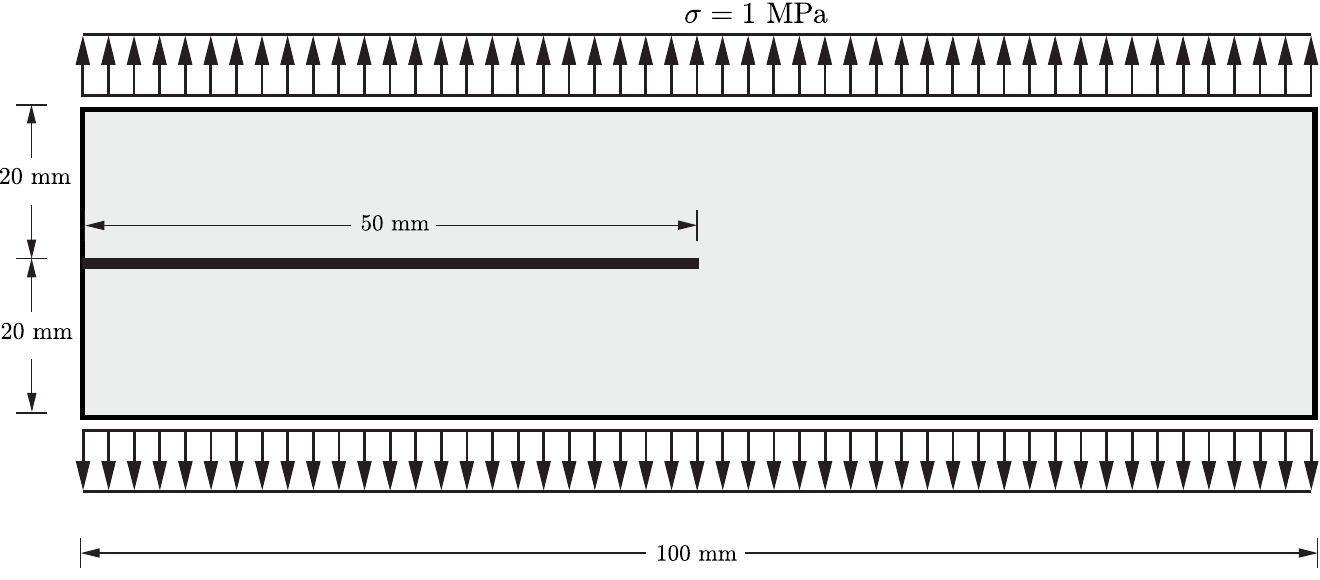}
    \caption{Dynamic crack branching. Dimensions and loading conditions.} 
    \label{fig:DynSketch}
\end{figure}

Finally, we conclude our study by examining the case of dynamic fracture, where inertia terms are present and off-diagonal matrices have a larger relative weight. The paradigmatic example of a rectangular specimen containing a sharp crack and subjected to a vertical tensile traction is considered, see Fig. \ref{fig:DynSketch} \cite{Borden2012}. The vertical traction is of magnitude $\sigma=1$ MPa and is instantaneously applied to the upper and lower boundaries. The dimensions of the specimen and the loading configuration are shown in Fig. \ref{fig:DynSketch}. The material parameters for the solid are set to $\rho=2450$ kg/$\tm{m}^3$, $E=32$ GPa, $\nu=0.2$, $\ell=0.25$ mm and $G_c=3$ J/$\tm{m}^2$, implying a Rayleigh wave speed of $v_r=2125$ m/s. The domain is uniformly discretised using linear quadrilateral elements with side length $h_e=0.25$ mm. 

\begin{figure}[H]
    \centering
    \includegraphics[width=\textwidth]{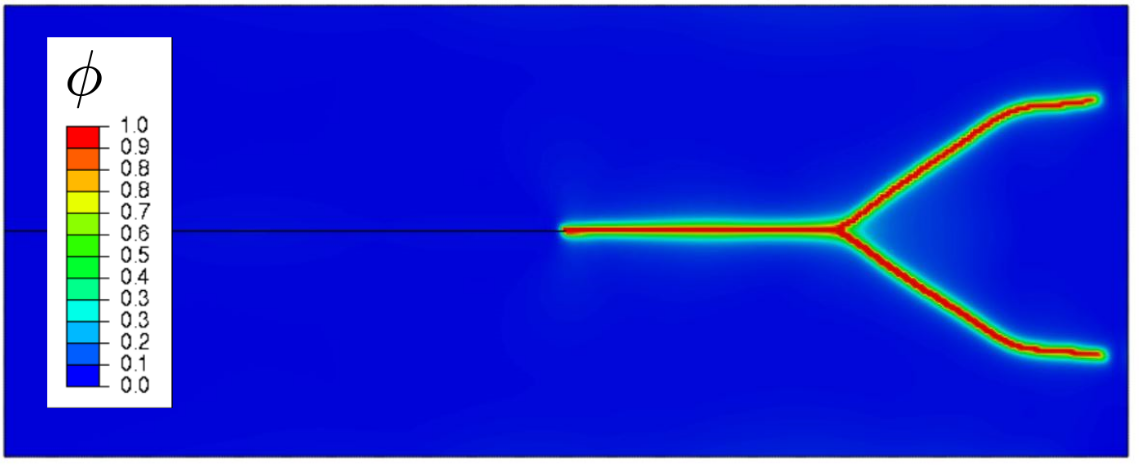}
    \caption{Dynamic crack branching. Crack trajectory predicted with the monolithic quasi-Newton implementation.}    
    \label{fig:DynCrack}
\end{figure}

The increment size for both the staggered and the monolithic calculations is set according to $\Delta t \approx h_e/v_r \approx 0.1\, \mu$s. Thus, differences in computation times are caused only by the computational cost of updating the Jacobian and from the number of iterations required for each increment to converge. The dynamic system is solved using a Backward Euler approach, without the need for using other algorithms such as HHT or Newmark's $\beta$-method to achieve convergence. The staggered and monolithic approaches show the same qualitative result, crack growth followed by branching - see Fig. \ref{fig:DynCrack}. However once the crack widening phase initiates, the staggered approach on average requires 4-6 times more iterations per increment, leading to a total computation time almost 3 times longer than the monolithic approach.

\begin{figure}[H]
    \centering
    \includegraphics[width=\textwidth]{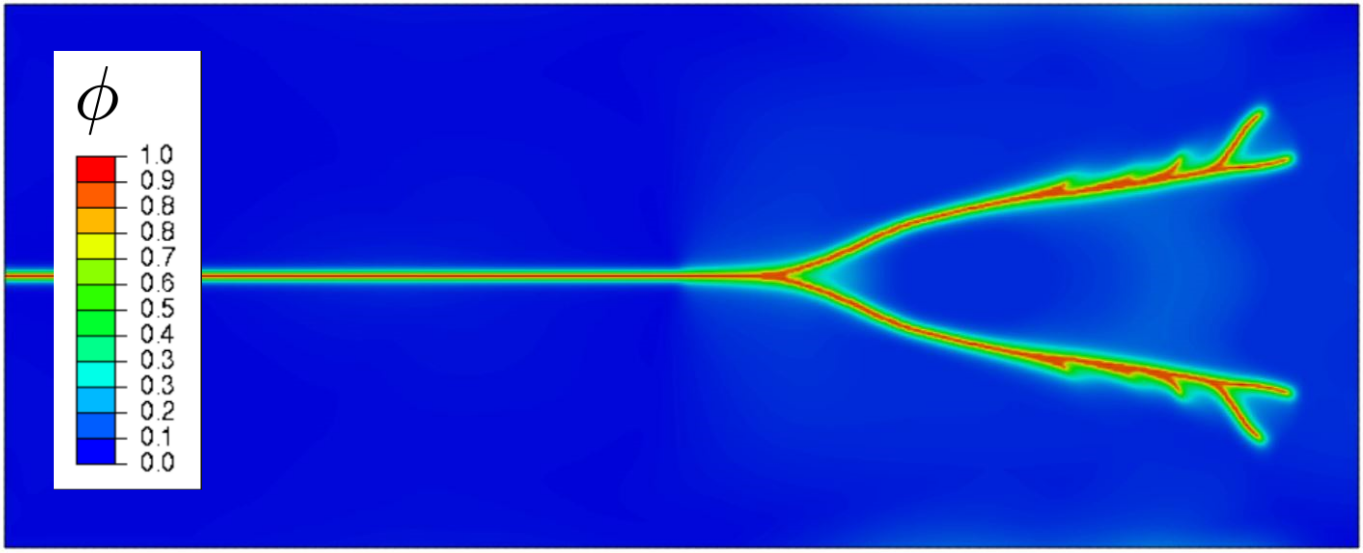}
    \caption{Dynamic crack branching. Crack trajectory predicted with the monolithic quasi-Newton implementation for $G_c=0.5$ J/$\tm{m}^2$.}    
    \label{fig:DynCrackGc}
\end{figure}

The result showcases how the monolithic quasi-Newton approach converges significantly faster than staggered solution schemes in highly non-linear problems. Moreover, it proves the capabilities of the quasi-Newton monolithic implementation in solving highly complex fracture problems. We further illustrate this aspect and the versatility of the method by obtaining results with $G_c=0.5$ J/$\tm{m}^2$, quadratic elements and an initial crack defined by prescribing $\phi$. As shown in Fig. \ref{fig:DynCrackGc}, complex crack patterns can be obtained with the present monolithic quasi-Newton implementation.

\section{Conclusions}
\label{Sec:Conclusions}

We present a quasi-Newton monolithic solution scheme for phase field fracture. The modelling framework makes use of the Broyden-Fletcher-Goldfarb-Shanno (BFGS) algorithm and is enhanced with a new adaptive time stepping algorithm. Several paradigmatic boundary value problems are solved, spanning stable and unstable quasi-static fracture, phase field fatigue and dynamic crack branching applications. Our main findings are:\\

\noindent (i) Monolithic quasi-Newton solution schemes are robust. Capable of solving various benchmark problems of varying complexity without convergence problems, unlike standard Newton monolithic frameworks.\\

\noindent (ii) By retaining unconditional stability, computation times are drastically reduced relative to widely used staggered solution schemes. Monolithic quasi-Newton computations are 10 to 100 times faster in all the problems considered. \\

\noindent (iii) Accurate phase field fatigue predictions can be obtained with only 4 increments per cycle, several orders of magnitude below the requirements of staggered approaches.\\

It is therefore expected that the use of monolithic quasi-Newton solution schemes will open new possibilities in phase field fracture modelling; for example, by enabling physically-based cycle-by-cycle fatigue calculations and tackling large scale problems. Of interest for future work is the study of the performance of the method in other challenging applications, including those involving large strains \cite{Reinoso2017a,Bilgen2019} and nearly incompressible materials \cite{Kumar2018}.

\section{Acknowledgments}
\label{Sec:Acknowledge of funding}

The authors gratefully acknowledge financial support from the Danish Hydrocarbon Research and Technology Centre (DHRTC) under the ``Reliable in-service assessment in aggressive environments'' project, publication
number DHRTC-PRP-108. E. Mart\'{\i}nez-Pa\~neda additionally acknowledges financial support from Wolfson College Cambridge (Junior Research Fellowship) and from the Royal Commission for the 1851 Exhibition through their Research Fellowship programme (RF496/2018).



\bibliographystyle{elsarticle-num} 
\bibliography{library}


\end{document}